\documentclass[10pt,psamsfonts]{amsart}
\usepackage{amsmath}
\usepackage{amsthm}
\usepackage{amssymb}
\usepackage{multicol}
\usepackage{amscd}
\usepackage{amsfonts}
\usepackage{amsbsy}
\usepackage{epsfig,afterpage}
\usepackage[dvips]{psfrag}
\usepackage{subfigure}
\usepackage{epsfig}
\usepackage{amsmath}
\usepackage[margin=3.6cm]{geometry}
\usepackage{srcltx}
\usepackage{amsfonts}
\usepackage{amssymb}
\usepackage{psfrag}
\usepackage{multicol}
\usepackage{verbatim}
\usepackage{graphicx}
\usepackage{amsmath}
\usepackage{amsthm}
\usepackage{amssymb}
\usepackage{amscd}
\usepackage{amsfonts}
\usepackage{amsbsy}
\usepackage{epsfig}
\usepackage[dvips]{psfrag}
\usepackage{multirow}
\usepackage{graphics}
\usepackage{epstopdf}
\usepackage{overpic}
\usepackage{float}
\usepackage[utf8]{inputenc}

\usepackage{graphicx}

\usepackage[colorlinks=true, linkcolor=blue, citecolor=red]{hyperref}

\def\e{\varepsilon}

\newcommand{\bbox}{\ \hfill\rule[-1mm]{2mm}{3.2mm}}

\newtheorem {theorem} {Theorem}%[section]
\newtheorem {proposition} [theorem]{Proposition}

\newtheorem {lemma}  [theorem]{Lemma}

\newtheorem{mtheorem}{Theorem}
\renewcommand*{\themtheorem}{\Alph{mtheorem}}

\usepackage{float}
\usepackage{tikz, wrapfig}
\tikzset{node distance=3cm, auto}
\allowdisplaybreaks

\begin{document}

\title[On the number of limit cycles for PWHS]
{On the number of limit cycles for piecewise polynomial holomorphic systems}

\author[Armengol Gasull, Gabriel Rondón and Paulo R. da Silva]{Armengol Gasull$^{1}$, Gabriel Rondón$^{2}$ and Paulo R. da Silva$^{2}$}
\address{$^{1}$Departament de Matemàtiques, Edifici Cc, Universitat Autònoma de Barcelona, 08193 Bellaterra, Barcelona, Spain; and
	Centre de Recerca Matem\`{a}tica, Campus de
	Bellaterra, 08193 Bellaterra, Barcelona, Spain}
\address{$^{2}$S\~{a}o Paulo State University (Unesp), Institute of Biosciences, Humanities and
	Exact Sciences. Rua C. Colombo, 2265, CEP 15054--000. S. J. Rio Preto, S\~ao Paulo,
	Brazil.}
 
\email{armengol.gasull@uab.cat}
\email{gabriel.rondon@unesp.br}
\email{paulo.r.silva@unesp.br}

\thanks{ .}

\makeatletter
\@namedef{subjclassname@2020}{\textup{2020} Mathematics Subject Classification}
\makeatother


\date{}
\maketitle

\begin{abstract}
In this paper we are concerned with determining lower bounds of the number of limit cycles for piecewise polynomial holomorphic systems with a straight line of discontinuity.
We approach this problem with  different points of view: study of the number of zeros of the first and second order averaging functions, or with the control of the limit cycles appearing from a monodromic equilibrium point via a degenerated Andronov--Hoph type bifurcation, adding at the very end the sliding effects. We also  use the Poincaré--Miranda theorem for obtaining an explicit piecewise linear holomorphic systems with 3 limit cycles, result that improves the known examples in the literature that had  a single limit cycle. 
\end{abstract}

%\tableofcontents

\section{Introduction}

The models in  nonsmooth dynamics of differential equations  have attracted the attention of many researchers  for the accuracy
of the obtained results comparing with the real
observations, see more details in the three books 
\cite{Acary, Brogliato, Kunze} and their references. Several of these models are given by piecewise smooth systems with some switching manifold. Moreover, on many of them the smooth systems are linear and the switching manifold is a straight line.

Holomorphic functions have a wide range of applications in several areas of applied science such as the study of fluid dynamics, for more information see, for instance, \cite{BatGK,Mars,Conw}. One of the most remarkable dynamical properties of holomorphic systems $\dot{z}=f(z)$ is the fact that these systems do not have limit cycles. The study and the properties of these systems make them interesting and beautiful but precisely this  absence of limit cycles makes them dynamically poor. However, in \cite{Rondon2022440} the authors proved that there are piecewise linear holomorphic systems that have one limit cycle. Moreover they proved that if the equilibrium points are on the straight line of discontinuity this limit cycle is unique.

In this paper we are interested  on piecewise polynomial holomorphic systems (PWHS). On one hand, each of the smooth systems has the beautiful properties of the holomorphic systems but on the other hand considered as a piecewise system they exhibits all the interesting features of the piecewise linear systems and much more. In particular, as we have already explained, they can have limit cycles.  Moreover, all the power of complex notation and of complex analysis can be used in their study.

Essentially, the techniques we employed to prove the results of this paper could be applied to piecewise smooth vector fields, without necessarily being piecewise holomorphic. However, the fact that it is holomorphic endows the vector field with important properties that simplify calculations. In fact, holomorphic systems, apart of the absence of limit cycles, have other surprising and interesting properties: reversibility, integrability, all their centers are isochronous, and complete knowledge of the phase portraits around their non-essential singularities, with simple local normal forms.  Moreover, when we add the hypothesis of being holomorphic to a system, we reduce the number of involved parameters. More concretely, while a polynomial system of degree $n$ depends on $n^2+3n+2$ parameters, a polynomial holomorphic system depends only on $2n+2$ parameters. In particular, while planar polynomial systems of degree $n$ that do not have continua of equilibrium  points have at most $n^2$ of them,  holomorphic polynomial  systems of degree $n$ have at most $n$ equilibrium points, provided that $f$ is not identically null.

The aim of this work is the study of the limit cycles of  piecewise polynomial holomorphic systems (PWHS), 
\begin{equation}\label{ch4:eq111}
\begin{aligned}
\left\{\begin{array}{l}
\dot{z}=F^{+}(z), \text{ when } \operatorname{Im}(z)> 0,\\[5pt]
\dot{z}=F^{-}(z),\text{ when } \operatorname{Im}(z)<0,
\end{array} \right.
\end{aligned}
\end{equation}
where $z=x+iy$ and $F^{\pm}(z)$ are holomorphic polynomial functions with $\operatorname{deg}(F^\pm)=n^\pm.$ 

Notice that the straight line $\Sigma=\{z\in\mathbb{C}:\operatorname{Im}(z)=0\}$ divides the plane in two half-planes 
$\Sigma^\pm$ given by $\{z\in\mathbb{C}:\operatorname{Im}(z)> 0\}$ and $\{z\in\mathbb{C}:\operatorname{Im}(z)< 0\}$, respectively. 
The orbits on $\Sigma$ are defined following the Filippov convention, see \cite{Filippov88}  for more details.

Let 
$\mathcal{L}_{n^+,n^-}\in\mathbb{N}
\cup\{\infty\}$ be the maximum number of limit cycles that \eqref{ch4:eq111} can have. 
The main goal of this article is to determine lower bounds for $\mathcal{L}_{n^+,n^-}.$ Notice, that trivially $\mathcal{L}_{n^+,n^-}=\mathcal{L}_{n^-,n^+}.$  The main tools that we will use are:

\begin{itemize}
\item The averaging theory up to second order to compute the so called {averaged functions of orders 1 and 2}. We will base our computations on the general results given in  \cite[Theorem 1]{MR3729598}. For the sake of completeness, we recall them  in Proposition \ref{prop3}. With this approach we obtain lower bounds of $\mathcal{L}_{n^-,n^+}$  (non optimal) for any $n^-,n^+.$ 

\item The computation and use of the Lyapunov quantities to produce limit cycle via degenerated Andronov--Hopf type bifurcations when $n^-+n^+\le 4.$  We  develop the method introduced in \cite{COLL2001671,1967ZaMM...47R.480R} and apply it in the particular holomorphic context to obtain explicit expression of the first five Lyapunov quantities.

\item The effect of the sliding for increasing by one the number of limit cycles obtained in the previous item. In this point we follow the ideas of \cite[Proposition 7.3]{doi:10.1137/11083928X} to perform a suitable final perturbation inside the holomorphic world.

\item  The use of the conformal map $w=1/z,$
to transform a neighborhood of infinity of our system with $n^+=n^-=1$ into a neighborhood of the origin of a new system \eqref{ch4:eq111} but with $n^+=n^-=2$. Then, by applying the results obtained by using the Lyapunov quantities in this latter case we prove that $\mathcal{L}_{1,1}\ge3$. This approach is very related with the one  of  Freire et al. \cite{MR4189023}, that proves the same result and with similar ideas but for piecewise linear systems.

\item The Poincaré--Miranda theorem to prove the existence of 3 limit cycles for an explicit piecewise linear holomorphic sytem, providing a second proof  that $\mathcal{L}_{1,1}\ge3$.  This result adapts to the holomorphic setting the same approach used in \cite{Armengol} to obtain a piecewise linear system with at least 3 limit cycles. 

\end{itemize}

Let us state our main results. Our first theorem deals with the analysis of the limit cycles obtained by computing  the second order averaged function of a perturbation of the global center $\dot z=iz$ inside system \eqref{ch4:eq111}. As usual, in next result, $[.]$ denotes the integer part of a real number, i.e. given $n\leq x<n+1$ then the integer part of $x$ is $[x]=n$.

\begin{mtheorem}\label{teo1}
Consider system \eqref{ch4:eq111} with $F^\pm(z)=iz+\epsilon h^{\pm}(z)$ and $\operatorname{deg}(h^\pm)=n^{\pm}.$  Then, by using the averaging theory up to order two, for $\epsilon$ small enough, the maximum number of limit cycles  that  bifurcate from the periodic orbits of the center is $\left[(3\max\{n^+,n^-\}-1))/2\right],$ this upper bound is attained and, in this case all these limit cycles are 
hyperbolic. As a consequence, $\mathcal{L}_{n^+,n^-}\ge \left[(3\max\{n^+,n^-\}-1))/2\right].$
\end{mtheorem}

Our second result studies with more detail the cases $n^++n^-\le 4$ and deals with the limit cycles obtained from a degenerated Andronov--Hopf bifurcation together with an extra limit cycle obtained by adding sliding (see Proposition \ref{prop_sliding}).  Of course this result is based on the computation of the Lyapunov quantities associated the weak focus-weak focus case, which are developed and  particularized for the piecewise holomorphic 
 systems.  We arrive until the 5th constant, see Theorem \ref{teob} for more details. The main difficulty in calculating these quantities arises from the length of the expressions involved, which was expected since in the discontinuous case we only consider half return maps around the monodromic equilibrium points and some cancellations that happen in the smooth case disappear. From this reason we have deferred some of the computations of our proof of Theorem \ref{teob} to an appendix. 

\begin{mtheorem}\label{teonew}
 Let $\mathcal{L}^0_{n^+,n^-}\le \mathcal{L}_{n^+,n^-}$ be the number of hyperbolic limit cycles bifurcating from the origin of system \eqref{ch4:eq111}. In Table \ref{tb_1} we give  some lower bounds of $\mathcal{L}^0_{n^+,n^-}$ for $n^++n^-\le 4,$
\begin{table}[h]
	\begin{center}
		\begin{tabular}{| c ||c| c | c |}
			\hline
			$n^+/n^-$& 1 & 2 & 3 \\
			\hline\hline
			1  & 1 & 3 & 5 \\
			\hline
			2  & 3 & 4 & - \\
			\hline
			3  & 5 & - & - \\
			\hline
		\end{tabular}
	\end{center}
	\vspace{0.2cm} \caption{Lower bounds for $\mathcal{L}^0_{n^+,n^-}.$ }\label{tb_1}
\end{table}
\end{mtheorem}

 We highlight that in the linear case $n^+=n^-=1,$ Theorems \ref{teo1} and \ref{teonew} only show that $\mathcal{L}_{1,1}\ge 1,$ that was a known result (\cite{Rondon2022440}). Our final theorem improves this lower bound.
 
\begin{mtheorem}\label{teod}
	There are piecewise linear holomorphic systems with at least 3 limit cycles, that is $\mathcal{L}_{1,1}\ge3.$
\end{mtheorem}

As we have already commented, we will present two proofs of the above theorem. The first one based on the study of a degenerate Andronov--Hopf bifurcation and the second one providing an explicit example and based on Poincaré--Miranda theorem. As we will see the first approach does not give explicit examples, but the existence of them. They exist when some involved parameters satisfy certain smallness properties, that  cannot be given explicitly.

As far as we know all examples with 3 limit cycles for piecewise linear systems have at least one of the involved smooth systems (or even both) of focus type. Since the real Jordan form of a focus can also be written as a holomorphic vector field, after all, it is not a full surprise that piecewise linear holomorphic systems do have also 3 limit cycles.

 More specifically, our second proof demonstrates that the piecewise linear holomorphic system 
       \begin{equation*}
\begin{aligned}
\left\{\begin{array}{l}
\dot{z}=\left(i+\frac{3}{8}\right)z-\frac{14333}{2000}+i\frac{1159}{1000}, \text{ when } \operatorname{Im}(z)> 0,\\[5pt]
\dot{z}=\left(i-\frac{1}{5}\right)z-\frac{51}{50}+i\frac{1}{250},\text{ when } \operatorname{Im}(z)<0,
\end{array} \right.
\end{aligned}
\end{equation*}
has 3 nested limit cycles that surround the real focus $(-1/5,-49/50).$ For more details see Section \ref{PM_3CL}. 

The construction the above example has been inspired by the example of piecewise linear system with 3 limit cycles presented by Llibre and Ponce in \cite{MR1673265} (developed by using Newton--Kantorovich theorem) and by the one of  Gasull and Mañosa in \cite{Armengol}  (also studied by using Poincaré--Miranda theorem). 

The paper is organized as follows. In Section \ref{sec:Preliminaries}, we present some basic results on the averaging theory, the Lyapunov method,  the Poincaré--Miranda theorem, and some more technical results. This section also contains Theorem \ref{teob}. Then we dedicate next three sections to prove Theorems \ref{teo1}, \ref{teonew}, and \ref{teod}. Finally, an appendix is also provided with some additional computations needed to prove Theorem \ref{teob}.

\section{Preliminaries}\label{sec:Preliminaries}
This section is devoted to establishing some results that will be used throughout the paper. We divide it in four subsections.

\subsection{The averaging method}\label{aver_metd}

We briefly recall some basic results of the averaging theory for piecewise smooth systems written in polar coordinates. An overview on this subject can be found in \cite{MR3729598}, and the reader can see the details of the proofs there. Consider the piecewise smooth systems of the form 
\begin{equation}\label{polar_sys}
\dfrac{dr}{d\theta}=\begin{aligned}
\left\{\begin{array}{l}
F^+(\theta,r,\epsilon) \text{ if } 0\leq\theta\leq\pi,\\[5pt]
F^-(\theta,r,\epsilon) \text{ if } \pi\leq\theta\leq2\pi,
\end{array} \right.
\end{aligned}
\end{equation}
where $F^\pm(\theta,r,\epsilon)=\sum_{i=1}^{k}\epsilon^iF_i^\pm(\theta,r)+\epsilon^{k+1}R^\pm(\theta,r,\epsilon),$ with $\theta\in S^1,$ $r>0$ and $\epsilon>0$ is a sufficiently small parameter. 

From \cite{MR3729598} we introduce the following functions $M_j^\pm(r)$ for $j=1,2$ related to the above system \eqref{polar_sys}:
$$M_j^\pm(r)=\dfrac{1}{j!}y_j(\pm\pi,r),\quad \text{for}\quad j=1,2,$$ where
\begin{itemize}
    \item $y_1^\pm(t,r)=\displaystyle\int_0^{t} F_1^\pm(\theta,r) d\theta$; 
    \item $y_2^\pm(t,r)=\displaystyle\int_0^{t} \left[2F_2^\pm(\theta,r)+2\partial_rF_1^\pm(\theta,r)y_1^\pm(\theta,r)\right] d\theta.$
\end{itemize}
Now, we define the function $M_j(r)=M_j^+(r)-M_j^-(r)$, which is called \textit{the averaged function of order $j.$} The following result can be found in \cite[Theorem 1]{MR3729598}:
\begin{proposition}\label{prop3} Let $M_l,$  $l\in\{1,2\},$ be the first non identically zero averaged function.  Then, each simple zero $r=r_0$ of  $M_{l}$ provides, for $\epsilon$ small enough, a hyperbolic limit cycle of the piecewise smooth system \eqref{polar_sys} that tends to $r=r_0$ when $\epsilon$ tends to 0.
\end{proposition}

    \subsection{The Lyapunov approach}\label{lya_metd}
    
To determine the existence of limit cycles associated with the system \eqref{ch4:eq111}, it is enough to calculate the zeros of the displacement function $\Delta_1(r)=(f^-)^{-1}(r)-f^+ (r )$, where $f^\pm$ are the half-return maps of $\dot z=F^\pm(z)$. 

Notice that these half return maps are locally well defined for instance when the equilibrium is a weak focus for both differential equations $\dot z= F^\pm(z)$ and it is located on  the line of discontinuity.

It is important to highlight that instead of the function $\Delta_1$ another function, that also captures the limit cycles,  can be defined using the composition of $f^+$ with $f^-$, that is, $\Delta_2(r)=f^-(f^+(r))-r.$ Both functions have the same zeros counting their multiplicities. Indeed, assume that $\rho^*$ is a zero of $\Delta_1$ of order $k\geq1$,  i.e. $\Delta_1^{(j)}(\rho^*)=0$ $j=0,\dots,k-1$ and $\Delta_1^{(k)}(\rho^*)\neq0.$ Recall that $\Delta_2$ can be rewritten as
$$\Delta_2(r)=f^-(f^+(r))-f^-(f^+(r)+\Delta_1(r)).$$
Thus, $\Delta_2^{(j)}(\rho^*)=0$ $j=0,\dots,k-1$ and $$\Delta_2^{(k)}(\rho^*)=-(f^-)'(f^+(\rho^*))\Delta_1^{(k)}(\rho^*)\neq0,$$ where we have used that $-(f^-)'(f^+(\rho^*))\neq0$ because both $f^\pm$ are diffeomorphisms. Therefore, $\rho^*$ is a zero of $\Delta_2$ of order $k.$ The reverse is done the same way.    
Both approaches are used throughout this work.    
    
Similarly, the main idea to determine the stability of an equilibrium point $p=(0,0)$ of equation \eqref{pwhs_eq_1} consists of starting from $(r_0,0)$, $r_0>0$ small enough, and to evaluate the sign of the first term nonzero of $\Delta_2$ in its power series expansion around $r=r_0.$ The coefficient of this term is usually called {\it a Lyapunov quantity or Lyapunov constant.}
    
 The Lyapunov quantities also give us information about the maximum number of limit cycles that bifurcate from an equilibrium  point of the PWHS, and this leads us to introduce the concept of order of a  weak focus. We say that $p=(0,0)$ is a \textit{weak focus of order $k$} of system \eqref{ch4:eq111}, whenever $F^\pm(0)=0$ and the displacement function defined as  $\Delta(r):=\Delta_2(r)=f^+(f^-(r))-r$ satisfies
 $$\Delta(r)=V_kr^k+\mathcal{O}(r^{k+1}), \quad\mbox{with}\quad V_k\ne0.$$   
As we have explained above, the coefficient $V_k$ is called 
 \textit{$k$th Lyapunov quantity}. In Proposition \ref{aux_lemma}, we will show that if the piecewise holomorphic system \eqref{ch4:eq111} has a weak focus order $k$ at the origin, then any small perturbation of system \eqref{ch4:eq111},  that keeps the singularity, has at most $k-1$ limit cycles. Moreover in that proposition we give effective computable conditions to check whether $k-1$ limit cycles do  appear. Finally, by introducing sliding in \eqref{ch4:eq111}, see Proposition \ref{prop_sliding}, we can obtain one extra limit cycle bifurcating from the origin. In short, from a weak focus of order $k,$ it is possible to bifurcate  $k$ limit cycles and there are verifiable conditions to ensure that this number of limit cycles do bifurcate from the origin.
 
In this context, just as in the smooth case it is possible to calculate the Lyapunov quantities, however in the discontinuous case there are also even Lyapunov quantities. The only difficulty in calculating these quantities arises from the length of the expressions involved. As we have already commented, this is because in the discontinuous case we only consider half return maps.     
    
Let us compute the first five Lyapunov quantities for piecewise holomorphic systems with a line of discontinuity following the ideas of \cite{COLL2001671,1967ZaMM...47R.480R}. Consider the system 
\begin{equation}\label{smooth_eq1}
\dot{z}=F(z)=\sum_{k=1}^{\infty}F_k(z),
\end{equation}
where $F_k(z)=A_k z^k$. Assume that the origin is of focus type and that the nearby solutions turn around it counterclockwise, that is, $\operatorname{Im}(a)>0,$ where $F_1(z)=az,$ $a\in\mathbb{C}.$

In the $(r,\theta)$-polar coordinates, system \eqref{smooth_eq1} is written as
\begin{equation*}
\begin{array}{rcl}
   \dfrac{dr}{d\theta}&=&r\dfrac{\cos(\theta)\operatorname{Re}(F(z)+\sin(\theta)\operatorname{Im}(F(z))}{\cos(\theta)\operatorname{Im}(F(z))-\sin(\theta)\operatorname{Re}(F(z))}\bigg|_{z=re^{i\theta}}\vspace{0.3cm}\\
   &=& r\dfrac{\operatorname{Re}(\overline{z}F(z))}{\operatorname{Im}(\overline{z}F(z))}\bigg|_{z=re^{i\theta}},
   \end{array}
\end{equation*}
or equivalently, 
\begin{equation}\label{smooth_eq_3}
    \dfrac{dr}{d\theta}=\dfrac{\sum_{k=1}^\infty r^k\operatorname{Re}(S_k(\theta))}{\sum_{k=1}^\infty r^{k-1}\operatorname{Im}(S_k(\theta))}=\sum_{k=1}^\infty R_k(\theta)r^k,
\end{equation}
where $S_k(\theta)=\overline{z}F(z)|_{z=e^{i\theta}}=e^{-i\theta}F_k(e^{i\theta})$ and the functions $R_k(\theta)$ can be easily obtained from them. Recall that $\operatorname{Im}(S_1(\theta))=\operatorname{Im}(a)>0.$

Let $r=R(\theta,s)$ be the solution of \eqref{smooth_eq_3} such that $R(0,s)=s,$ and write it  as
\begin{equation*}\label{polar_sol_1}
    R(\theta,s)-s=\sum_{k=1}^\infty \omega_k(\theta)s^k,\quad\text{where}\quad\omega_k(0)=0,\quad\text{for}\quad k\geq 1.
\end{equation*}
Plugging  this expression in \eqref{smooth_eq_3} we get successive linear differential equations for the unknowns $\omega_k$ and we can obtain them.
Then, in the smooth case, the Lyapunov quantities are given by the first $k$ such that $\omega_k(2\pi)\ne0$. It is known that this $k$ is odd, see for instance \cite{Andronov}.

In the sequel we will perform a similar computation but for the piecewise holomorphic system
 \begin{equation}\label{pwhs_eq_1}
\begin{aligned}
\left\{\begin{array}{l}
\dot{z}=F^{+}(z), \text{ when } \operatorname{Im}(z)> 0,\\[5pt]
\dot{z}=F^{-}(z),\text{ when } \operatorname{Im}(z)<0,
\end{array} \right.
\end{aligned}
\end{equation}
where $F^\pm(z)=\sum_{k=1}^\infty F_k^\pm(z)$ are defined by \eqref{smooth_eq1} and $F_1^\pm(z)=(i+\lambda^\pm)z.$

As we will see, the essential difference with the smooth case is that the Lyapunov quantities can be obtained by the values $\omega_k(\pi),$ computed for both smooth differential systems $\dot z=F^\pm(z).$ 

 \begin{mtheorem}\label{teob}
     Consider system \eqref{pwhs_eq_1} where \
     $F^\pm(z)=\sum_{k=1}^{\infty}F_k^\pm(z).$ Suppose that $F_1^\pm(z)=(i+\lambda^\pm)z$ and write  $F_2^\pm(z)=A^\pm z^2,$ $F_3^\pm(z)=B^\pm z^3,$ $F_4^\pm(z)=C^\pm z^4,$ and $F_5^\pm(z)=D^\pm z^5.$ Then its first five Lyapunov quantities are:
     \begin{itemize}
         \item[(i)] $V_1=e^{(\lambda^++\lambda^-)\pi}-1;$ 
         \item[(ii)] $V_2=\omega_2^+(\pi)+\omega_2^-(\pi)e^{3\lambda^+\pi};\vspace{0.3cm}$ 
         \item[(iii)] $V_3=e^{\lambda^+\pi}\omega_3^+(\pi)-2(\omega_2^+(\pi))^2+\omega_3^-(\pi)e^{5\lambda^+\pi};\vspace{0.3cm}$
          \item[(iv)] $V_4=e^{2\lambda^+\pi}\omega_4^+(\pi)-5e^{\lambda^+\pi}\omega_2^+(\pi)\omega_3^+(\pi)+5(\omega_2^+(\pi))^3+\omega_4^-(\pi)e^{7\lambda^+\pi};\vspace{0.3cm}$
           \item[(v)] $\!\!\!\begin{array}{rcl}
           V_5\!\!\!\!&=&\!\!\!\!e^{3\lambda^+\pi}\omega_5^+(\pi)+21e^{\lambda^+\pi}(\omega_2^+(\pi))^2\omega_3^+(\pi)-14(\omega_2^+(\pi))^4\vspace{0.3cm}\\
           & &\!\!\!\!-3e^{2\lambda^+\pi}(\omega_3^+(\pi))^2-6e^{2\lambda^+\pi}\omega_2^+(\pi)\omega_4^+(\pi)+\omega_5^-(\pi)e^{9\lambda^+\pi}.
           \end{array}$
     \end{itemize}
     where $\omega_i^\pm(\pi)$ for $i=1,\cdots,5$ are:

      $$\begin{array}{rcl}
      \omega_2^\pm(\pi)&=&e^{\lambda^\pm\pi}(-e^{\lambda^\pm\pi}-1)\operatorname{Re}\left[ \pm\dfrac{(1+\lambda^\pm i)A^\pm}{i+\lambda^\pm}\right];\vspace{0.3cm}\\%%%
    \omega_3^\pm(\pi)&=&e^{\lambda^\pm\pi}(e^{2\lambda^\pm\pi}-1)\left\{\operatorname{Re}\left[\dfrac{(1+\lambda^\pm i)B^\pm}{2(i+\lambda^\pm)}+\dfrac{ A^\pm\overline{A^\pm}}{4}\right]\right.\vspace{0.2cm}\\
     & &- \left.\operatorname{Im}\left[\dfrac{(1+\lambda^\pm i)(A^\pm)^2}{4(i+\lambda^\pm)}\right]\right\}+e^{-\lambda^\pm\pi}(\omega_2^\pm(\pi))^2;\vspace{0.3cm}\\%%%%%
     \omega_4^\pm(\pi)
&=&-2e^{-2\lambda^\pm\pi}\left(\omega_2^\pm(\pi)\right)^3+3e^{-\lambda^\pm\pi}\omega_2^\pm(\pi)\omega_3^\pm(\pi)\vspace{0.2cm}\\
&&\mp\dfrac{1}{2}e^{\lambda^\pm\pi}(-e^{3\pi\lambda^\pm}-1)\left\{\operatorname{Re}\left[\eta_1^\pm-\lambda^\pm\eta_5^\pm-2(\gamma_1^\pm+\gamma_2^\pm)\right]\right.\vspace{0.2cm}\\
 &&-\left.\operatorname{Im}\left[\lambda^\pm\eta_1^\pm+\eta_3^\pm+2(\gamma_3^\pm+\gamma_4^\pm))\right]\right\}\vspace{0.2cm}\\
 &&\pm\dfrac{1}{2}e^{\lambda^\pm\pi}(-e^{\pi\lambda^\pm}-1)\left\{\operatorname{Re}\left[\eta_2^\pm-\lambda^\pm\eta_6^\pm\right]-\operatorname{Im}\left[\lambda^\pm\eta_2^\pm+\eta_4^\pm\right]\right\};\vspace{0.3cm}\\%%%%%%%%%%%
 \end{array}$$
  $$\begin{array}{rcl}
     \omega_5^\pm(\pi)
&=&-\dfrac{5}{2}e^{-3\lambda^\pm\pi}(\omega_2^\pm(\pi))^4+2e^{-2\lambda^\pm\pi}(\omega_2^\pm(\pi))^2\omega_3^\pm(\pi)+\dfrac{3}{2}e^{-\lambda^\pm\pi}(\omega_3^\pm(\pi))^2\vspace{0.2cm}\\
 &&+e^{\lambda\pi}(e^{4\pi\lambda^\pm}-1)\left\{\operatorname{Re}\left[\xi_1^\pm+\xi_3^\pm+\xi_5^\pm+\xi_8^\pm+\xi_9^\pm+\xi_{10}^\pm+\xi_{11}^\pm+\xi_{12}^\pm-2(\delta_1^\pm+\delta_3^\pm)\right]\right.\vspace{0.2cm}\\
&&+\left.\operatorname{Im}\left[\xi_6^\pm+\xi_7^\pm-2(\delta_6^\pm-\lambda^\pm\delta_1^\pm)\right]\right\}\pm e^{\lambda\pi}(-e^{4\pi\lambda^\pm}-1)\left\{\operatorname{Re}\left[\xi_2^\pm\right]-\operatorname{Im}\left[\xi_4^\pm\right]\right\}\vspace{0.2cm}\\
 &&\mp\omega_2^\pm(\pi)(-e^{3\pi\lambda^\pm}-1)\left\{\operatorname{Re}\left[\eta_1^\pm-\lambda^\pm\eta_5^\pm-4(\gamma_1^\pm+\gamma_2^\pm)\right]\right.\vspace{0.2cm}\\
 &&-\left.\operatorname{Im}\left[\lambda^\pm\eta_1^\pm+\eta_3^\pm+4(\gamma_3^\pm+\gamma_4^\pm)\right]\right\}\vspace{0.2cm}\\
&& +\dfrac{1}{4}e^{\lambda\pi}(e^{2\pi\lambda^\pm}-1)\left\{\operatorname{Re}\left[B^\pm\delta_7^\pm+\dfrac{\lambda^\pm}{2}\delta_{10}^\pm\right]-\dfrac{1}{2}\operatorname{Im}\left[(A^\pm)^2\delta_7^\pm\right]\right\}\vspace{0.2cm}\\
&&-\dfrac{1}{4}e^{\lambda\pi}(-e^{\pi\lambda^\pm}-1)\left\{\operatorname{Re}\left[B^\pm\delta_8^\pm+\dfrac{\lambda^\pm}{2}\delta_9^\pm-8(\delta_2^\pm+\delta_4^\pm)\right]\right.\vspace{0.2cm}\\
&&-\left.\dfrac{1}{2}\operatorname{Im}\left[(A^\pm)^2\delta_8^\pm+16(\delta_5^\pm-\lambda^\pm\delta_2^\pm)\right]\right\}\\
&&\pm (-e^{\pi\lambda^\pm}-1)\omega_2^\pm(\pi)\left\{\operatorname{Re}\left[\eta_2^\pm-\lambda^\pm\eta_6^\pm\right]-\operatorname{Im}\left[\lambda^\pm\eta_2^\pm+\eta_4^\pm\right]\right\};
      \end{array}$$
for certain constants $\eta^\pm_{k\geq 1},\gamma^\pm_{k\geq 1},\xi^\pm_{k\geq 1},$ and $\delta^\pm_{k\geq 1}$ given in the appendix. Moreover, the stability of origin of system \eqref{smooth_eq1} is determined by the sign of the first non zero Lyapunov quantity.   
     \end{mtheorem}

To prove the previous theorem we will use the following result, whose proof can be consulted in \cite{MR1460166}. In the sequel, we use the notation $\widetilde{f}=\widetilde{f}(\theta)=\int_{0}^\theta f(s)ds.$ 
    \begin{proposition}\label{coef_R} 
    	Consider the analytic differential equation
    	\begin{equation}\label{smooth_eq_7_}
    		\dfrac{dr}{d\theta}=\sum_{k=2}^\infty R_k(\theta)r^k.
    	\end{equation}
     Let $r=R(\theta,\rho)$ be its solution with initial condition $R(0,\rho)=\rho.$ Then 
    	\begin{equation}\label{polar_sol_1}
    		R(\theta,\rho)-\rho=\sum_{k=2}^\infty u_k(\theta)\rho^k,
    	\end{equation}
    where  $u_k(0)=0,$  $k\geq 2,$ and
        $$\begin{array}{rcl}
            u_2(\theta)&=&\widetilde{R}_2; \\
            u_3(\theta)&=&(\widetilde{R}_2)^2+\widetilde{R}_3;\\
            u_4(\theta)&=&(\widetilde{R}_2)^3+2\widetilde{R}_2\widetilde{R}_3+\widetilde{\widetilde{R}_2R_3}+\widetilde{R}_4=(\widetilde{R}_2)^3+3\widetilde{R}_2\widetilde{R}_3-\widetilde{\widetilde{R}_3R_2}+\widetilde{R}_4;\\
            u_5(\theta)&=&(\widetilde{R}_2)^4+3(\widetilde{R}_2)^2\widetilde{R}_3+\widetilde{(\widetilde{R}_2)^2R_3}+2\widetilde{R}_2\widetilde{\widetilde{R}_2R_3}+\frac{3}{2}(\widetilde{R}_3)^2+2\widetilde{R}_2\widetilde{R}_4+2\widetilde{R_4\widetilde{R}_2}+\widetilde{R}_5\\
            &=&(\widetilde{R}_2)^4+5(\widetilde{R}_2)^2\widetilde{R}_3+\widetilde{(\widetilde{R}_2)^2R_3}-2\widetilde{R}_2\widetilde{\widetilde{R}_3R_2}+\frac{3}{2}(\widetilde{R}_3)^2+4\widetilde{R}_2\widetilde{R}_4-2\widetilde{\widetilde{R}_4R_2}+\widetilde{R}_5.\\
        \end{array}$$
    \end{proposition}
The next proposition is a technical result, which is proven in the Appendix.
     \begin{proposition}\label{main:prop}
     Let $r=R(\theta,s)$ be the solution with initial condition  $R(0,s)=s$  of system $
     	\dot{z}=F(z)=\sum_{k=1}^{\infty}F_k(z).$ Assume that $F_1(z)=(\lambda+i)z, \lambda\in\mathbb{R},$ and write $F_2(z)=Az^2,$ $F_3(z)=Bz^3,$ $F_4(z)=Cz^4,$ and $F_5(z)=Dz^5.$ Then 
     $$R(\pi,s)-s=\omega_1(\pi)s+\omega_2(\pi)s^2+\omega_3(\pi)s^3+\omega_4(\pi)s^4+\omega_5(\pi)s^5+\mathcal{O}(s^6),$$
     where 
     $$\begin{array}{rcl}
    \omega_1(\pi)&=&e^{\lambda\pi}-1;\vspace{0.3cm}\\
     \omega_2(\pi)&=&e^{\lambda\pi}(-e^{\lambda\pi}-1)\operatorname{Re}\left[ \dfrac{(1+\lambda i)A}{i+\lambda}\right]\vspace{0.3cm}; \\ 
    \omega_3(\pi)&=&e^{\lambda\pi}(e^{2\lambda\pi}-1)\left\{\operatorname{Re}\left[\dfrac{(1+\lambda i)B}{2(i+\lambda)}+\dfrac{ A\overline{A}}{4}\right]-\operatorname{Im}\left[\dfrac{(1+\lambda i)A^2}{4(i+\lambda)}\right]\right\}+e^{-\lambda\pi}(\omega_2(\pi))^2;\vspace{0.3cm}\\
     \omega_4(\pi)
&=&-2e^{-2\lambda\pi}\left(\omega_2(\pi)\right)^3+3e^{-\lambda\pi}\omega_2(\pi)\omega_3(\pi)\\
&&-\dfrac{1}{2}e^{\lambda\pi}(-e^{3\pi\lambda}-1)\left\{\operatorname{Re}\left[\eta_1-\lambda\eta_5-2(\gamma_1+\gamma_2)\right]-\operatorname{Im}\left[\lambda\eta_1+\eta_3+2(\gamma_3+\gamma_4)\right]\right\}\vspace{0.3cm}\\
 &&+\dfrac{1}{2}e^{\lambda\pi}(-e^{\pi\lambda}-1)\left\{\operatorname{Re}\left[\eta_2-\lambda\eta_6\right]-\operatorname{Im}\left[\lambda\eta_2+\eta_4\right]\right\};\vspace{0.3cm}\\
%         \end{array}$$
%     $$\begin{array}{rcl}
\omega_5(\pi)
&=&-\dfrac{5}{2}e^{-3\lambda\pi}(\omega_2(\pi))^4+2e^{-2\lambda\pi}(\omega_2(\pi))^2\omega_3(\pi)+\dfrac{3}{2}e^{-\lambda\pi}(\omega_3(\pi))^2\vspace{0.3cm}\\
 &&+e^{\lambda\pi}(e^{4\pi\lambda}-1)\left\{\operatorname{Re}\left[\xi_1+\xi_3+\xi_5+\xi_8+\xi_9+\xi_{10}+\xi_{11}+\xi_{12}-2(\delta_1+\delta_3)\right]\right.\vspace{0.3cm}\\
&&+\left.\operatorname{Im}\left[\xi_6+\xi_7-2(\delta_6-\lambda\delta_1
)\right]\right\}+e^{\lambda\pi}(-e^{4\pi\lambda}-1)\left\{\operatorname{Re}\left[\xi_2\right]-\operatorname{Im}\left[\xi_4\right]\right\}\vspace{0.3cm}\\
    
 &&-\omega_2(\pi)(-e^{3\pi\lambda}-1)\left\{\operatorname{Re}\left[\eta_1-\lambda\eta_5-4(\gamma_1+\gamma_2)\right]-\operatorname{Im}\left[\lambda\eta_1+\eta_3+4(\gamma_3+\gamma_4)\right]\right\}\vspace{0.3cm}\\
&& +\dfrac{1}{4}e^{\lambda\pi}(e^{2\pi\lambda}-1)\left\{\operatorname{Re}\left[B\delta_7+\dfrac{\lambda}{2}\delta_{10}\right]-\dfrac{1}{2}\operatorname{Im}\left[A^2\delta_7\right]\right\}\vspace{0.3cm}\\
&&-\dfrac{1}{4}e^{\lambda\pi}(-e^{\pi\lambda}-1)\left\{\operatorname{Re}\left[B\delta_8+\dfrac{\lambda}{2}\delta_9-8(\delta_2+\delta_4)\right]-\dfrac{1}{2}\operatorname{Im}\left[A^2\delta_8+16(\delta_5-\lambda\delta_2)\right]\right\}\vspace{0.3cm}\\
&&+\omega_2(\pi)(-e^{\pi\lambda}-1)\left\{\operatorname{Re}\left[\eta_2-\lambda\eta_6\right]-\operatorname{Im}\left[\lambda\eta_2+\eta_4\right]\right\};\\
      \end{array}$$
      for certain constants $\eta_{k\geq 1},\gamma_{k\geq 1},\xi_{k\geq 1},$ and $\delta_{k\geq 1}$ given in the appendix.
     \end{proposition}
     \begin{proof}[Proof of Theorem \ref{teob}]
To prove our result we have to compose the map induced by the flow of $\dot{z}=F^+(z)$ between $\theta=0$ and $\theta=\pi$, namely, $f^+$ and the map induced by the flow of $\dot{z}=F^-(z)$ between $\theta=\pi$ and $\theta=2\pi$, namely, $f^-$.    From Proposition \ref{main:prop} by substituting $F$ by $F^+$, we know that
$$\begin{array}{rcl}
f^+(s)&=&(\omega_1^++1)s+\omega_2^+s^2+\omega_3^+s^3+\omega_4^+s^4+\omega_5^+s^5+\mathcal{O}(s^6)\vspace{0.3cm}\\
&:=&f_1s+f_2s^2+f_3s^3+f_4s^4+f_5s^5+\mathcal{O}(s^6).
\end{array}$$
     The second one can be computed also by Proposition \ref{main:prop}. Thus, we get that
     $$\begin{array}{rcl}
     f^-(s)&=&(\omega_1^-+1)s+\omega_2^-s^2+\omega_3^-s^3+\omega_4^-s^4+\omega_5^-s^5+\mathcal{O}(s^6)\vspace{0.3cm}\\
     &:=&g_1s+g_2s^2+g_3s^3+g_4s^4+g_5s^5+\mathcal{O}(s^6),
     \end{array}$$
where $\omega_k^-$ are given in the proposition substituting $F$ by $-F^-(-z).$ Hence, the stability of the origin is controlled by 
$$\begin{array}{rcl}
\Delta(s)&:=& f^-(f^+(s))-s\vspace{0.3cm}\\
&=&(g_1f_1-1)s+(g_1f_2+g_2f_1^2)s^2+(g_1f^3+2g_2f_1f_2+g_3f_1^3)s^3\vspace{0.3cm}\\
& &+(g_4f_1^4+3g_3f_1^2f_2+g_2f_2^2+2g_2f_1f_3+g_1f_4)s^4\vspace{0.3cm}\\
& &+(g_5f_1^5+4g_4f_1^3f_2+3g_3f_1f_2^2+3g_3f_1^2f_3+2g_2f_2f_3+2g_2f_1f_4+g_1f_5)s^5\vspace{0.3cm}\\
& &+\mathcal{O}(s^6)
\end{array}
$$ and the Lyapunov quantities of \eqref{pwhs_eq_1} are
$$\begin{array}{rcl}
V_1&=&g_1f_1-1;\vspace{0.3cm}\\
V_2&=&g_1f_2+g_2f_1^2,\,\text{when}\, V_1=0;\vspace{0.3cm}\\
V_3&=&g_1f^3+2g_2f_1f_2+g_3f_1^3,\,\text{when}\, V_1=V_2=0;\vspace{0.3cm}\\
V_4&=&g_4f_1^4+3g_3f_1^2f_2+g_2f_2^2+2g_2f_1f_3+g_1f_4,\,\text{when}\, V_1=V_2=V_3=0;\vspace{0.3cm}\\
V_5&=&g_5f_1^5+4g_4f_1^3f_2+3g_3f_1f_2^2+3g_3f_1^2f_3+2g_2f_2f_3+2g_2f_1f_4+g_1f_5,\, V_i=0, i=1,2,3,4.
\end{array}$$     
Then by substituting in the above formulas the values obtained from Proposition \ref{main:prop} the proof follow.
     \end{proof}
 
    We end this section by presenting our version of a  well-known method to obtain the maximum number of limit cycles born from a weak-focus of order $k$ of system \eqref{pwhs_eq_1} via the so called {\it degenerated Andronov--Hopf bifurcation}.  Of course it will use the expressions of the Lyapunov quantities given  in Theorem \ref{teob}.
  
 Let $s_j,$ $j=1,\cdots,k-1$ be $k-1$ real parameters and consider the piecewise perturbed holomorphic systems:
 \begin{equation}\label{pwhs_eq_per_1}
 	\begin{aligned}
 		\left\{\begin{array}{l}
 			\dot{z}=G^+(z,s), \text{ when } \operatorname{Im}(z)> 0,\\[5pt]
 			\dot{z}=G^-(z,s),\text{ when } \operatorname{Im}(z)<0,
 		\end{array} \right.
 	\end{aligned}
 \end{equation}
 where $G^\pm(z,\mathbf{0})=F^\pm(z),$ $s=(s_1,\cdots,s_{k-1})\in\mathbb{R}^{k-1}$  and $\mathbf{0}\in\mathbb{R}^{k-1}$. We denote ${W_j}(s)$ the Lyapunov quantities associated with system \eqref{pwhs_eq_per_1}. Recall that ${W_j}(\mathbf{0})=V_j,$ for all $j=1,\cdots,k-1$. As usual, given a differentiable map $g:\mathbb{R}^n\to\mathbb{R}^n$ we denote by $J(g)$ its Jacobian matrix.
 
 \begin{proposition}\label{aux_lemma} 
 	If piecewise holomorphic system  \eqref{pwhs_eq_per_1}, with $s=0$ has a weak focus order $k$ at $p=(0,0),$ then at most $k-1$ limit cycles (counting their multiplicities) bifurcate from $p.$ 
 	Moreover, if $\det(J_{s}({W_1},\cdots,{W_{k-1}}))(0)\neq 0$, then there exist real parameters $s_j,$ $j=1,\cdots,k-1$ small enough such it  has $k-1$ hyperbolic limit cycles bifurcating from the origin.   
 \end{proposition}
 \begin{proof}
 	Assume without loss of generality that for system  \eqref{pwhs_eq_per_1}, with $s=0,$ $p=(0,0)$ is a repelling focus of order $k$, i.e. $V_k>0.$
 	
 	The proof of the first part of the theorem is by contradiction. Thus, suppose that there exist $\rho_j(s)>0$ for $j=1,\cdots,k$ such that $D(\rho_j(s),s)=0,$ where $D(r,s)=\Delta(r,s)/r$ being $\Delta(r,s)$ the displacement function associated with the piecewise perturbed holomorphic system \eqref{pwhs_eq_per_1}. Notice that $\rho_j(s)\rightarrow 0,$ when $s\rightarrow 0.$ Moreover, 
 	we can suppose that $\rho_j(s)\neq \rho_l(s)$ for all $l,j=1,\cdots,k.$
 	
 	Consider the functions $g_j(r,s)=\dfrac{\partial^jD}{\partial r^j}(r,s),$ for $j=1,\cdots,k-1.$ From Rolle's theorem, for each $j,$ there exist $\mu_{j,l}(s)>0,$ for $l=1,\cdots,k-j$ such that $g_{j}(\mu_{j,l}(s),s)=0$ and $\lim_{s\to 0}\mu_{j,l}(s)=0.$ For $j=k-1,$ $g_{k-1}(\mu(s),s)=0,$ where $\mu(s)=\mu_{k-1,1}(s)$ and
 	$$g_{k-1}(r,s)=(k-1)!W_{k}(s)+\mathcal{O}(r).$$
 	In particular, 
 	$$0=\lim_{s\to 0}g_{k-1}(\mu(s),s)=(k-1)!V_{k}\ne0,$$
  giving the desired contradiction.
 	
 	Now, we prove the second part of the theorem. Since ${W_k}(\mathbf{0})=V_k>0$, then there exists $r_k$ such that $0<r_{k}\ll1$ and $\Delta(r_{k},\mathbf{0})>0$.
 	
 	Consider the  smooth function $f:\mathbb{R}^{k-1}\rightarrow\mathbb{R}^{k-1}$ defined by 
 	$$f(s)=({W_1}(s),\cdots,{W_{k-1}}(s)),$$ where $f(\mathbf{0})=\mathbf{0}$. Since $\det(J_{s}f)(0)\neq 0$, then from \textit{Inverse Function Theorem} there exist neighborhoods $\mathcal{U}\subset\mathbb{R}^{k-1}$ of $s=\mathbf{0}$ and  $\mathcal{V}\subset\mathbb{R}^{k-1}$ of $f(\mathbf{0})=\mathbf{0}$ such that $f(\mathcal{U})\subset\mathcal{V}$ and $f:\mathcal{U}\rightarrow\mathcal{V}$ is bijective. Thus, there exist real parameters $s_j^*,$ $j=1,\cdots,k-1$ small enough such that ${W_{k-1}}(s^*)<0,$ where $s^*=(s_{1}^*,\cdots,s_{k-1}^*).$ In addition, by continuity ${W_k}(s^*)={W_k}(\mathbf{0})+\mathcal{O}(s^*)>0$ and $\Delta(r_{k},s^*)>0.$ Since the displacement function associated with piecewise perturbed holomorphic system \eqref{pwhs_eq_per_1} is given by
 	$$\Delta(r,s^*)=W_{k-1}(s^*)r^{k-1}+\mathcal{O}(r^{k}),\quad\mbox{with}\quad W_{k-1}(s^*)<0,$$
 	then there exists $r_{k-1}$ such that $r_{k-1}\ll r_k$ and  $\Delta(r_{k-1},s^*)<0$.
 	From the Intermediate Value Theorem, there exists $\rho^*_{k-1}\in(r_{k-1},r_{k})$ such that $\Delta(\rho^*_{k-1},s^{*})=0.$ Therefore, system \eqref{pwhs_eq_per_1} has at least one limit cycle.
 		
 	As before, using the bijectivity of $f$, we know that there exist real parameters $s_j^{**},$ $j=1,\cdots,k-1$ small enough such that ${W_{k-2}}(s^{**})>0,$ where $s^{**}=(s_{1}^{**},\cdots,s_{k-1}^{**}).$ Moreover, by continuity ${W_k}(s^{**})>0,$ $\Delta(r_{k},s^{**})>0,$ ${W_{k-1}}(s^{**})<0$ and $\Delta(r_{k-1},s^{**})<0.$ Since the displacement function associated with piecewise perturbed holomorphic system \eqref{pwhs_eq_per_1} is given by
 	$$\Delta(r,s^{**})=W_{k-2}(s^{**})r^{k-2}+\mathcal{O}(r^{k-1}), \quad\mbox{with}\quad W_{k-2}(s^{**})>0,$$
 	then there exists $r_{k-2}$ such that $r_{k-2}\ll r_{k-1}$ and  $\Delta(r_{k-2},s^{**})>0$. Again, from the Intermediate Value Theorem, there exist $\rho^{**}_{k-2}\in(r_{k-2},r_{k-1})$ and $\rho^{**}_{k-1}\in(r_{k-1},r_{k})$ such that $\Delta(\rho^{**}_{k-2},s^{**})=0$ and $\Delta(\rho^{**}_{k-1},s^{**})=0.$ Therefore, system \eqref{pwhs_eq_per_1} has at least two limit cycles.
 	
 	We can do this process $k-1$ times, thus the piecewise perturbed holomorphic system \eqref{pwhs_eq_per_1} has $k-1$ limit cycles bifurcating from the origin, which must be hyperbolic because in the first part of the proposition we have proved that $k-1$ is the sum of the multiplicities of all the limit cycles bifurcating from the origin.
 \end{proof}
 
   \subsection{Poincaré--Miranda theorem}
  The Poincaré--Miranda theorem was formulated and proven by H. Poincaré in 1883 and 1886 respectively, \cite{P1,P2} and this arises as an extension of the Bolzano theorem to higher dimensions. In 1940, C. Miranda obtained the same result as an equivalent formulation of Brouwer's fixed point theorem, \cite{Miranda}.
\begin{theorem}[Poincaré--Miranda]\label{PM}
Consider the set
$$B=\{x=(x_1,\cdots, x_n)\in\mathbb{R}^n: L_i<x_i<U_i, 1\leq i\leq n\}.$$ Suppose that $f=(f_1,f_2,\cdots,f_n):\overline{B}\subset\mathbb{R}^n\to\mathbb{R}^n$ is continuous, $f(x)\neq 0$ for all $x\in\partial B$, and for $1\leq i \leq n$,
$$f_i(x_1,\cdots,x_{i-1},L_i,x_{i+1},\cdots,x_n)\leq 0\quad \text{and}\quad f_i(x_1,\cdots, x_{i-1},U_i, x_{i+1},\cdots,x_n)\geq 0.$$
Then, there exists $s\in B$ such that $f(s)=0$.
\end{theorem} 

\begin{center}
	\begin{figure}[h]
		\begin{overpic}[width=9truecm,height=3truecm]{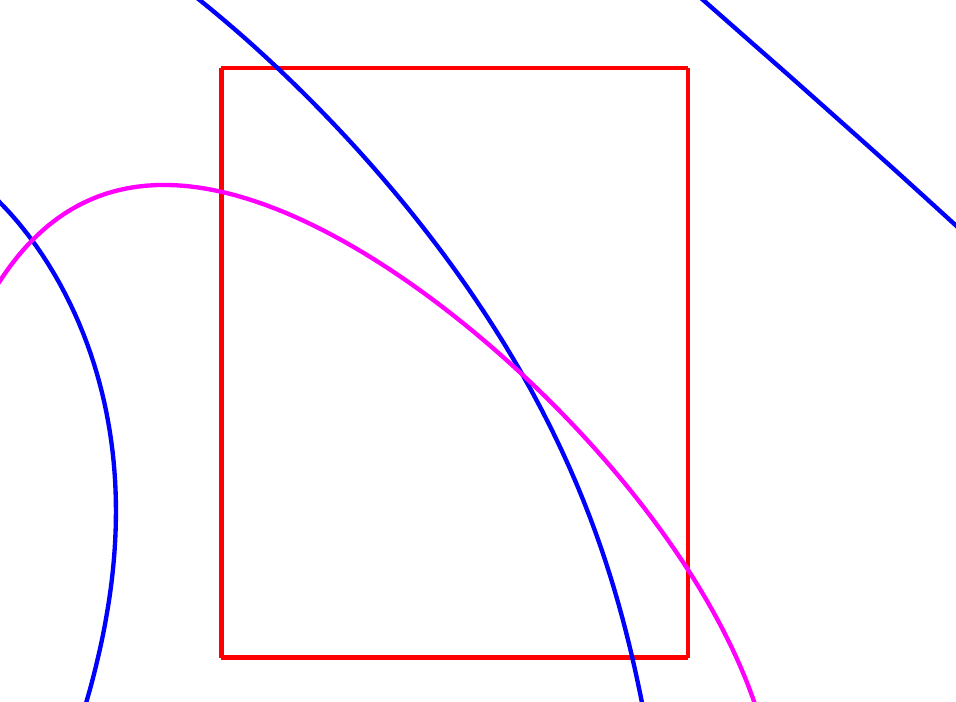}
			\put(1,30){{\color{blue} $f_1(x,y)=0$}}
			\put(77,4){{\color{magenta}$f_2(x,y)=0$}}
			\put(11.5,15){{\color{red}$f_1<0$}}
			\put(73,15){{\color{red} $f_1>0$}}
			\put(45,32){{\color{red}$f_2>0$}}
			\put(45,-2){{\color{red} $f_2<0$}}
		\end{overpic}
		\caption{Illustration of a box $B$  where Poincaré--Miranda theorem applies for $n=2.$}
	\end{figure}
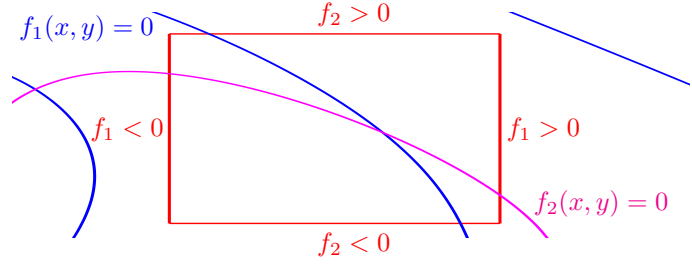
\end{center}

\subsection{A miscellany of results}

This first lemma will be useful to prove the existence of simple zeros of the averaged functions. Its proof can be found in \cite[Lemma 4.5]{MR2170413}.

\begin{lemma}\label{lemma1}
	Consider $p+1$ linearly independent functions $f_i:U\subset\mathbb{R}\rightarrow\mathbb{R},$ $i=0,\cdots,p$
	\begin{itemize}
		\item[(a)] Given $p$ arbitrary values $x_i\in U,$ $i=0,\cdots,p$ there exist $p+1$ constants $C_i,$ $i=0,\cdots,p$ such that
		\begin{equation}\label{eq:lemma}
			f(x):=\sum_{i=0}^p C_if_i(x)
		\end{equation}
		is not the zero function and $f(x_i)=0,$ $i=0,\cdots,p.$
		\item[(b)] Furthermore, if all $f_i$ are analytic functions on $U$ and there exist $j=0,\cdots,p$ such that $f_j|_U$ has a constant sign, it is possible to get an $f$ given by \eqref{eq:lemma}, such that it has at least $p$ simple zeros in $U.$
	\end{itemize}
\end{lemma}
 
Another result that we will use in this work is the celebrated  \textit{Descartes Theorem}, which gives us information about the number of positive zeros of a real polynomial in terms of the sign and number of monomials, for more details see for instance \cite{MR0174165}. Given an ordered list of $p+1$ non-zero real numbers $[a_0,a_1,\ldots,a_p]$ we will say that its number of sign variations  is $m,$ $0\le m\le p,$ if there are exactly $m$ values of $j\le p-1$ such that $a_j a_{j+1}<0.$ 

\begin{theorem}[Descartes Theorem]\label{descartes}
	Consider the real polynomial $P(x)=a_{0}x^{i_0}+\cdots+a_{p}x^{i_p}$ with $0\leq i_0<\cdots<i_p$ and  $a_{j}$ non-zero real constants for $j\in\{0,\ldots,p\}.$ If the number of sign variations  of $[a_0,a_1,\ldots,a_p]$  is $m$, then $P(x)$ has exactly   $m-2n$ positive real zeros counting their multiplicities, where $n$ is a non negative integer number.
\end{theorem}

To finish this section we state the following lemma, as it was established in \cite{Armengol}. It is a consequence of Taylor's formula. We will use it to control the signs of the functions appearing in the faces of the boxes when we apply Poincaré--Miranda theorem.
 
\begin{lemma}\label{lemma2}
Consider the function $h(x)=A\cos(\alpha x)+B\sin(\alpha x) +Ce^{\beta x}+De^{-\beta x}$, with $A,B,C,D\in\mathbb{R}$, $\alpha\neq0$, $\beta> 0$ and $x\in[\underline{x},\overline{x}]\subset\mathbb{R}^+$. Then for each $n\geq0$ we have $h(x)=\displaystyle\sum_{j=0}^n a_jx_j + m_n(x)x^{n+1}$, where
$$a_j =\frac{1}{j!}\left(\alpha_j\left[A\cos\left(j\frac{\pi}{2}\right)+ B \sin\left(j\frac{\pi}{2}\right)\right]\right) + \beta_j\left[C+(-1)^j D\right],$$
\begin{equation}\label{M}
    |m_n(x)|\leq\overline{m}=\frac{|\alpha|^{n+1}(|A|+|B|)+|\beta|^{n+1}(|C|e^{\beta \overline{x}} + |D|e^{-\beta \underline{x}})}{(n + 1)!}.
\end{equation}
\end{lemma}
\section{Proof of Theorem \ref{teo1}}\label{sec:mainresults}

Consider PWHS \eqref{ch4:eq111} such that $F^\pm(z)=iz+\epsilon h^{\pm}(z)$, where $h^{\pm}(z)=\sum_{k=0}^{n^\pm}(a^\pm_k+ib^\pm_k)z^k,$ for some $n^\pm\geq 1$ and $z=x+iy.$ In the $(r,\theta)-$polar coordinates $z=re^{i\theta}$, \eqref{ch4:eq111} is converted into 
\begin{equation}\label{polar_eq1}
\begin{aligned}
\left\{\begin{array}{l}
\dfrac{dr}{d\theta}=\epsilon\dfrac{a_1^+ r^2+r c_0^{+}+\sum_{k=1}^{n^+-1}r^{k+2}c_k^+}{r(1+b_1^+ \epsilon)+\epsilon d_0^{+}+\epsilon\sum_{k=1}^{n^+-1}r^{k+1}d_k^\pm}=F^+(r,\theta,\epsilon), \text{ in } \Omega_{0,\pi},\vspace{0.5cm}\\ [5pt]
\dfrac{dr}{d\theta}=\epsilon\dfrac{a_1^- r^2+r c_0^{-}+\sum_{k=1}^{n^--1}r^{k+2}c_k^-}{r(1+b_1^- \epsilon)+\epsilon d_0^{-}+\epsilon\sum_{k=1}^{n^--1}r^{k+1}d_k^-}=F^-(r,\theta,\epsilon),\text{ in } \Omega_{\pi,2\pi},
\end{array} \right.
\end{aligned}
\end{equation}
where 
$$\begin{array}{ll}
\Omega_{0,\pi}=\{(r,\theta):\,r\geq 0\, \text{and}\,\theta\in[0,\pi]\};&d_0^\pm=b_{0}^\pm\cos(\theta)-a_{0}^\pm\sin(\theta);\\
\Omega_{\pi,2\pi}=\{(r,\theta):\,r\geq 0\, \text{and}\,\theta\in[\pi,2\pi]\};&c_k^\pm=a_{k+1}^\pm\cos(k\theta)-b_{k+1}^\pm\sin(k\theta);\\ 
c_0^\pm=a_{0}^\pm\cos(\theta)+b_{0}^\pm\sin(\theta);&d_k^\pm=b_{k+1}^\pm\cos(k\theta)+a_{k+1}^\pm\sin(k\theta).\\
\end{array}$$
for all $ k\geq 1.$ Thus, expanding $F^\pm$ around $\epsilon=0,$ \eqref{polar_eq1} is written as
\begin{equation*}\label{polar_eq11}
    \dfrac{dr}{d\theta}=\epsilon F_1^\pm(\theta,r)+\epsilon^2F_2^\pm(\theta,r)+\mathcal{O}(\epsilon^3),
\end{equation*}
where 
$$
\begin{array}{rcl}
 F_1^\pm(\theta,r)&=&a_1^\pm r+c_0^\pm+\displaystyle\sum_{k=1}^{n^\pm-1}r^{k+1}c_k^\pm;\vspace{0.3cm}\\
 F_2^\pm(\theta,r)&=&\dfrac{1}{r}\left(-b_1^\pm r-d_0^\pm-\displaystyle\sum_{k=2}^n d_{k-1}^\pm r^k\right)\left(a_1^\pm+c_0^\pm+\displaystyle\sum_{k=2}^n c_{k-1}^\pm r^k\right).
\end{array}$$
Computing the first averaged function
\begin{equation*}
   M_1^\pm(r)=\displaystyle\int_0^{\pm\pi} F_1^\pm(\theta,r) d\theta= 2b_0^\pm \pm a_1^\pm r\pi-2\displaystyle\sum_{k=1}^{[\frac{n^\pm}{2}]}\frac{b_{2k}^\pm}{2k-1}r^{2k}.
\end{equation*}
Now, without loss of generality suppose that $n^+\geq n^-,$ then
\begin{equation*}
\begin{array}{rcl}
    M_1(r)&=&M_1^+(r)-M_1^-(r)\\
        &=&2(b_0^+-b_0^-)+\pi(a_1^{+}+a_1^-
)r+2\displaystyle\sum_{k=1}^{[\frac{n^-}{2}]}\frac{(b_{2k}^{-}-b_{2k}^{+})}{2k-1}r^{2k}-2\displaystyle\sum_{k=[\frac{n^-}{2}]+1}^{[\frac{n^+}{2}]}\frac{b_{2k}^{+}}{2k-1}r^{2k}.
\end{array}
\end{equation*}

The function $M_1(r)$ is a polynomial with  $\left[n^+/2\right]+2$ monomials. Then, by Descartes' Theorem (see Theorem \ref{descartes}),  $M_1(r)$ has at most $\left[n^+/2\right]+1$ positive zeros, taking into account their multiplicities. Moreover, since the coefficients of the monomials can be chosen independently we conclude that the coefficients of $M_1(r)$ are totally free. Finally, we can 
use Lemma \ref{lemma1}, because any the functions $r^{2k}$ have a constant sign. As a consequence we can choose the parameters of the perturbation in such a way that $M_1(r)$ has exactly   $\left[n^+/2\right]+1$ simple zeros. Hence, from Proposition \ref{prop3},  the same number of hyperbolic limit cycle bifurcate from the periodic orbits of $\dot z=iz,$ for $|\epsilon|$ small enough.

To get more limit cycles we need to go through the computation of the second averaged function. For that, we need that the averaged function of first order be identically zero, thus $a_1^+=-a_1^-,$ $b^+_{2k}=b^-_{2k},$ for all $k=0,\cdots,\left[{n^-}/{2}\right]$ and $b^+_{2k}=0$ for all $k=\left[{n^-}/{2}\right]+1,\cdots,\left[{n^+}/{2}\right].$ Then,  

\begin{equation*}
\begin{array}{rcl}
   M_2^\pm(r)&=&\dfrac{1}{2}\displaystyle\int_0^{\pm\pi} \left[2F_2^\pm(\theta,r)+2\partial_rF_1^\pm(\theta,r)y_1^\pm(\theta,r)\right] d\theta\vspace{0.3cm}\\
            &=& 4a_1^\pm a_0^\pm-2b_0^\pm b_1^\pm \pm\pi(-a_1^\pm b_1^\pm-2b_0^\pm a_2^\pm-2b_2^\pm a_0^\pm)r\vspace{0.3cm}\\
            &&+\displaystyle\sum_{k=1}^{n^\pm-1}U_k^{\pm}r^{2k}\pm\pi a_1^\pm\displaystyle\sum_{k=1}^{\left[\frac{n^\pm-1}{2}\right]}b_{2k+1}^\pm r^{2k+1},\vspace{0.3cm}
\end{array}
\end{equation*}
with
\begin{align*}
	U_k^\pm:=-\frac{2^{k+1} k!}{(2k)!}\sum_{2s+t=2k+1}\left(\dfrac{(2(k-1))!}{2^{k-2}(k-1)!}a_{2s}^\pm a_t^\pm+\dfrac{(t-2s)(2k-2)!}{(2s-1)2^{k-1}(k-1)!}b_{2s}^\pm b_t^\pm\right).
\end{align*}
Thus, 
$$
\begin{array}{rcl}
      M_2(r)&=& M_2^+(r)-M_2^-(r)\vspace{0.2cm}\\
            &=& -4a_1^-(a_0^++a_0^-)-2b_0^-(b_1^+-b_1^-)+\pi\left(a_1^-(b_1^+-b_1^-)\right.\vspace{0.2cm}\\
            &&-\left.2b_0^-(a_2^++a_2^-)-2b_2^-(a_0^++a_0^-)\right)r+\displaystyle\sum_{k=1}^{n^--1}(U_k^+-U_k^-)r^{2k}+\displaystyle\sum_{k=n^-}^{n^+-1}U_k^+ r^{2k}\vspace{0.2cm}\\
            &&+\pi a_1^-\displaystyle\sum_{k=1}^{\left[\frac{n^--1}{2}\right]}(b_{2k+1}^--b_{2k+1}^+) r^{2k+1}-\pi a_1^-\displaystyle\sum_{k=\left[\frac{n^-+1}{2}\right]}^{\left[\frac{n^+-1}{2}\right]}b_{2k+1}^+ r^{2k+1}.
            \end{array}
$$
%with,
%\begin{align*}
%	S_k:=-\frac{2^{k+1} k!}{(2k)!}&\sum_{2s+t=2k+1}\left(\dfrac{(2(k-1))!}{2^{k-2}(k-1)!}(a_{2s}^+ a_t^+-a_{2s}^- a_t^-)\right. \vspace{0.2cm}\\
%	&\quad+\left.\dfrac{(t-2s)(2k-2)!}{(2s-1)2^{k-1}(k-1)!}(b_{2s}^+ b_t^+-b_{2s}^- b_t^-)\right),
% \end{align*}
%\begin{align*}
%	T_k:=-\frac{2^{k+1} k!}{(2k)!} \sum_{2s+t=2k+1}\left(\dfrac{(2(k-1))!}{2^{k-2}(k-1)!}a_{2s}^+ a_t^+\dfrac{(t-2s)(2k-2)!}{(2s-1)2^{k-1}(k-1)!}b_{2s}^+ b_t^+\right).
%\end{align*}

In this case $M_2(r)$ is a polynomial with
exactly $\left[{(3\max\{n^+,n^-\}+1)}/{2}\right]$ monomials. Moreover, again its coefficients can be taken with total freedom by choosing in a suitable way the coefficients of the perturbation. Hence we can follow the same steps that in the study of $M_1(r)$  and prove that $M_2(r)$ can have
exactly $\left[{(3\max\{n^+,n^-\}-1)}/{2}\right]$  simple zeros.  Then, applying  Proposition \ref{prop3} the result follows.

\section{Proof of Theorem \ref{teonew}}\label{sec:tB}

As we will see, our proof of this theorem will be a straightforward consequence of the four propositions proved in this section. 

 \begin{proposition}\label{imp_prop3}
       Consider the PWHS
       \begin{equation}\label{pw_fam1}
\begin{aligned}
\left\{\begin{array}{l}
\dot{z}=(i+\lambda+s_1)z+(A^{+}+s_2)z^2, \text{ when } \operatorname{Im}(z)> 0,\\[5pt]
\dot{z}=(i-\lambda)z,\text{ when } \operatorname{Im}(z)<0,
\end{array} \right.
\end{aligned}
\end{equation}
where $\lambda\neq 0$, $A^+=a_1+ib_1$, $a_1=\dfrac{b_1(\lambda^2-1)}{\lambda}$ and $b_1\neq 0$. Then there exist real parameters $s_1$ and $s_2$ small enough such that it has two limit cycles bifurcating from the origin.
   \end{proposition}
   \begin{proof}
       Computing the Lyapunov functions and expanding them around $s_1=s_2=0,$ we get the following expressions, where $s=(s_1,s_2)$:
       \begin{itemize}
           \item $W_1=\pi s_1+\mathcal{O}_2(s)$; \vspace{0.2cm}
           \item $W_2=-\dfrac{2e^{\pi\lambda}(e^{\pi\lambda}+1) \lambda}{1+\lambda^2}s_2+\dfrac{b_1e^{\pi\lambda}(e^{\pi\lambda}+1)}{\lambda}s_1+\mathcal{O}_2(s)$;\vspace{0.2cm}
           \item $W_3=\dfrac{b_1^2e^{2\pi\lambda}(e^{2\pi\lambda}-1)(1+\lambda^2)}{8\lambda^2}+\mathcal{O}_1(s)$; \vspace{0.2cm}
       \end{itemize}
       Notice that $W_1=W_2=0$ and $W_3\neq 0$ provided that $s=0.$ Moreover, 
       $$\det(J_s(W_1,W_2))(0)=-\dfrac{2\pi e^{\pi\lambda}(e^{\pi\lambda}+1) \lambda}{1+\lambda^2}\neq 0.$$
       From Proposition \ref{aux_lemma}, we can conclude that there exist real parameters $s_1$ and $s_2$ small enough such that system \eqref{pw_fam1} has 2 limit cycles.
   \end{proof}
   \begin{proposition}\label{imp_prop2}
       Consider the PWHS
       \begin{equation}\label{pw_fam2}
\begin{aligned}
\left\{\begin{array}{l}
\dot{z}=(i+\lambda+s_1)z+(s_2+is_3)z^2+(B^++s_4)z^3, \text{ when } \operatorname{Im}(z)> 0,\\[5pt]
\dot{z}=(i-\lambda)z,\text{ when } \operatorname{Im}(z)<0,
\end{array} \right.
\end{aligned}
\end{equation}
where $\lambda\neq 0,$ $B^+=a_2+ib_2,$ $a_2=\dfrac{b_2(-1+\lambda^2)}{2\lambda}$ and $b_2\neq 0$. Then there exist real parameters $s_1,s_2,s_3$ and $s_4$ small enough such that it has four limit cycles bifurcating from the origin.
   \end{proposition}
    \begin{proof}
       Computing the Lyapunov functions and expanding them around $s_1=s_2=s_3=s_4=0,$ we get the following expressions, where $s=(s_1,s_2,s_3,s_4)$:
       \begin{itemize}
           \item $W_1=\pi s_1+\mathcal{O}_2(s)$;\vspace{0.2cm}
           \item $W_2=\dfrac{e^{\pi\lambda} (e^{\pi\lambda}+1) (-1+\lambda^2)}{1+\lambda^2}s_3-\dfrac{2e^{\pi\lambda} (e^{\pi\lambda}+1)\lambda}{1 + \lambda^2}s_2+\mathcal{O}_2(s)$;\vspace{0.2cm}
           \item $W_3=\dfrac{e^{2\pi\lambda} (e^{2\pi\lambda}-1)\lambda}{1+\lambda^2}s_4-\dfrac{
  b_2 e^{2\pi\lambda} (e^{2\pi\lambda}-1)}{2\lambda}s_1+\mathcal{O}_2(s)$;\vspace{0.2cm}
           \item $W_4=-\dfrac{b_2 e^{3\pi\lambda} (  e^{3\pi\lambda}+1) (7+27 \lambda^2)}{12(1+9\lambda^2)}s_3+\dfrac{b_2e^{3\pi\lambda}( e^{3\pi\lambda}+1)(1+21\lambda^2)}{12\lambda(1+9\lambda^2)}s_2+\mathcal{O}_2(s)$;\vspace{0.2cm}
           \item $W_5=-\dfrac{b_2^2 e^{4\pi\lambda} (e^{4\pi\lambda}-1)(\lambda^2-1)}{2(1 + \lambda^2)}+\mathcal{O}_1(s).$\vspace{0.2cm}
       \end{itemize}
       Notice that $W_1=W_2=W_3=W_4=0$ and $W_5\neq 0$ provided that $s=0.$ Moreover, 
       $$\det(J_s(W_1,W_2,W_3,W_4))(0)=\dfrac{b_2 e^{6\pi\lambda} (e^{\pi\lambda}-1)(e^{\pi\lambda}+1)^3(1-e^{\pi\lambda} + e^{2\pi\lambda})\pi] (1+33\lambda^2)}{12(1+\lambda^2)(1+9\lambda^2)}\neq 0.$$
         From Proposition \ref{aux_lemma}, we can conclude that there exist real parameters $s_1,s_2,s_3$ and $s_4$ small enough such that system \eqref{pw_fam2} has 4 limit cycles.
   \end{proof}
    \begin{proposition}\label{imp_prop}
       Consider the PWHS
       \begin{equation}\label{pw_fam3}
\begin{aligned}
\left\{\begin{array}{l}
\dot{z}=(i+1+s_1)z+(A^++s_2)z^2, \text{ when } \operatorname{Im}(z)> 0,\\[5pt]
\dot{z}=(i-1)z+(A^-+is_3)z^2,\text{ when } \operatorname{Im}(z)<0,
\end{array} \right.
\end{aligned}
\end{equation}
where 
\begin{multicols}{2}
\begin{itemize}
    \item $A^\pm=\mp1+ib_1^\pm;$
    \item $b^\pm_1=\dfrac{p^*\mp\tau(\tau\pm2)}{2\tau}$;
    \item $\tau\in\mathbb{R}\setminus\{0\};$
    \item $p^*=4+8\displaystyle\sum_{k=1}^5e^{k\pi}+4e^{6\pi}.$
\end{itemize}
\end{multicols}
\noindent Then there exist real parameters $s_1,s_2$ and $s_3$ small enough such that it has 3 limit cycles bifurcating from the origin.
   \end{proposition}
   \begin{proof}
    Computing the Lyapunov functions and expanding them around $s_1=s_2=s_3=0,$ we get the following expressions, where $s=(s_1,s_2,s_3)$:
       \begin{itemize}
           \item $W_1=\pi s_1+\mathcal{O}_2(s)\vspace{0.3cm}$; 
           \item $
            W_2=\dfrac{e^\pi}{2\tau}\left((1+e^\pi)(p^*-\tau(2+2\pi+\tau))-2\pi \tau\right)s_1\vspace{0.3cm}-e^\pi(1+e^\pi)s_2+\mathcal{O}_2(s);
           $
           \item $
           W_3=\dfrac{1}{16\tau^2}P_1(\tau,4)s_1+\dfrac{e^{2\pi}}{16\tau^2}P_2(\tau,3)s_2-\dfrac{e^{2\pi} (e^{2\pi}-1)(\tau^2+p^*)}{4\tau}s_3+\mathcal{O}_2(s);$\vspace{0.3cm}
           \item $
           W_4=\dfrac{e^{-12\pi}(1+e^\pi)}{120\tau}P_3(\tau,4)+\mathcal{O}_1(s).
          $\vspace{0.3cm}
       \end{itemize}
   where $P_j(\tau,n),$ $j=1,2,3,$ are polynomials in $\tau$ of degree $n$ with coefficients in   $\mathbb{Q}[e^\pi],$ with very large expressions and we do not explicit them. 
        
       It is easy to see that $P_3(\tau,4)$ does  not have real zeros.  Thus, $W_4\neq 0$ for all $\tau\in\mathbb{R}\setminus\{0\}.$ Moreover, $W_1=W_2=W_3=0$ and $W_4\neq 0$ provided that $s=0$ and $\tau\in\mathbb{R}\setminus\{0\}.$ Finally, 
       $$
       \det(J_s(W_1,W_2,W_3))(0)=
       \dfrac{\pi e^{3\pi}}{4\tau}(e^{\pi}-1)(e^{\pi}+1)^2(\tau^2+p^*)\neq 0,$$
      for all $\tau\in\mathbb{R}\setminus\{0\}.$ 
      
       Hence, from Proposition \ref{aux_lemma}, we can conclude that there exist real parameters $s_1,s_2$ and $s_3$ small enough such that system \eqref{pw_fam3} has 3 limit cycles.
   \end{proof}

   \noindent\textbf{Remark.}
       Note that in Proposition \ref{imp_prop} the Lyapunov quantities depend on a parameter $\tau.$ This is because such quantities involve polynomials and square roots of expressions which make more difficult to study the order of focus of system \eqref{pw_fam3}. For that reason, we have used rational parameterizations to facilitate such a study. For more details see, for instance, \cite{GasullTorre}.

       Before obtaining the coefficients of system \eqref{pw_fam3} that lead us to obtain the focus of order 4, we calculate the Lyapunov quantity $V_3$ of system \eqref{pw_fam3} for $r=s=q=0, $ we get
       $$V_3=\dfrac{1}{4}e^{2\pi}(e^{\pi}-1)\left((b_1^+)^2+2a^-_1 b_1^++p^*(a^-_1)^2-2a^-_1b_1^--(b_1^-)^2\right).$$
       To determine the order of the focus we must solve $V_3=0,$ which leads us to obtain 
       $$b_1^-=-a^-_1\pm\sqrt{(b_1^++a^-_1)^2+p^*(a^-_1)^2}$$
       Determining the values of $V_4$ with this root leads us to ugly  expressions. An alternative way to solve $V_3=0$ is to seek for a new variable $\tau$. To find such change of variables we set 
       $$M^2=L^2+p^*(a^-_1)^2,\quad \text{with}\quad L=b_1^++a^-_1.$$
       Since $M^2-L^2=(M+L)(M-L)$, thus we introduce a real parameter $\tau$ such that
       $$M+L=\tau\quad\text{and}\quad M-L=\dfrac{p^*(a^-_1)^2}{\tau}.$$
       Solving for $L$ and $M$ we get that
       $$M=\dfrac{p^*(a^-_1)^2+\tau^2}{2\tau}\quad\text{and}\quad L=\dfrac{p^*(a^-_1)^2-\tau^2}{2\tau}.$$
       This implies that $V_3$ is zero provided that 
       $$b_1^-=-a_1^-\pm \dfrac{p^*(a^-_1)^2+\tau^2}{2\tau}.$$
       Consequently, all the coefficients are rational expressions that depend on the new parameter $\tau.$ From them we obtain that $W_4$  is also a rational expression involving $\tau.$ This computational trick allows to prove in a much more easy and computable way the existence of this weak focus of order 4.

Observe that in the above three propositions $F^\pm(0)=0$, as in  system \eqref{pwhs_eq_1}, and hence  we have not considered sliding segments. In next proposition, by introducing a new real parameter $d$ small enough such that $F^-(0)\ne0,$ we get a sliding segment and it is possible to obtain one more limit cycle with this mechanism, which also bifurcates from the origin. To prove the proposition we adapt to the piecewise holomorphic setting the ideas of the approach done in \cite{doi:10.1137/11083928X} where the authors study the first return map near the origin for  piecewise linear differential equations.

\begin{proposition}\label{prop_sliding}
Consider the piecewise perturbed holomorphic system 
       \begin{equation}\label{pw_fam4}
\begin{aligned}
\left\{\begin{array}{l}
\dot{z}=(i+\lambda^+)z+\displaystyle\sum_{k=2}^\infty A^+_kz^k, \text{ when } \operatorname{Im}(z)> 0,\\[5pt]
\dot{z}=(i+\lambda^+)d+(i+\lambda^-)z+\displaystyle\sum_{k=2}^\infty A^-_kz^k,\text{ when } \operatorname{Im}(z)<0,
\end{array} \right.
\end{aligned}
\end{equation}
where $d\in\mathbb{R}$ is a small parameter. If $\lambda^++\lambda^->0$ (resp. $\lambda^++\lambda^-<0$), then there exists $d<0$ (resp. $d>0$) sufficiently small such that one limit cycle bifurcates from the origin.
 \end{proposition}
\begin{proof}
For $d=0$, the first Lyapunov quantity is given by $$V_1=e^{\pi(\lambda^++\lambda^-)}-1.$$
    Since $\lambda^+\neq-\lambda^-,$ then the origin is an attractor or repulsive equilibrium point. Hence, using the same ideas of \cite[Proposition 7.3]{doi:10.1137/11083928X}, we get that the displacement function for $d$ small enough is given by
$$\Delta(r,d)=P(r)-r=d(e^{\lambda^-\pi}+1)+W_1(d)r+\mathcal{O}(r^2),$$
where $W_1(0)=V_1$ and $P$ is the first return map associated with system \eqref{pw_fam4}. Since $d$ is a sufficiently small arbitrary parameter, we can use the same ideas of the proof of Proposition \ref{aux_lemma} and conclude that a limit cycle bifurcates from the origin.
\end{proof}

\begin{proof}[Proof of Theorem \ref{teonew}]
   
By using the results proved in   Propositions \ref{imp_prop3}, \ref{imp_prop2} and \ref{imp_prop} we have the results shown in Table \ref{tb_2}.

\begin{table}[h]
	\begin{center}
		\begin{tabular}{| c ||c| c | c |}
			\hline
			$n^+/n^-$& 1 & 2 & 3 \\
			\hline\hline
			1  & 0 & 2 & 4 \\
			\hline
			2  & 2 & 3 & - \\
			\hline
			3  & 4 & - & - \\
			\hline
		\end{tabular}
	\end{center}
	\vspace{0.2cm} \caption{Lower bounds for $\mathcal{L}^0_{n^+,n^-}$ without sliding.}\label{tb_2}
\end{table}

Finally,  by perturbing the piecewise holomorphic systems \eqref{pw_fam1}, \eqref{pw_fam2} and \eqref{pw_fam3}, given in these propositions,  as in Proposition \ref{prop_sliding}, we obtain that these perturbed systems have an additional limit cycle bifurcating fom the origin. Consequently, the number of limit cycles found for these systems are summarized in the Table \ref{tb_1}, given in the statement of the theorem. \end{proof}

\section{Proof of Theorem \ref{teod}}\label{sec:examples}

 \subsection{First proof of Theorem \ref{teod}}\label{example2}
    Consider the piecewise linear holomorphic system
       \begin{equation}\label{pw_3cycles}
\begin{aligned}
\left\{\begin{array}{l}
\dot{z}=-(i-1)z-A^--is_3, \text{ when } \operatorname{Im}(z)> 0,\\[5pt]
\dot{z}=-(i+1+s_1)z-A^+-s_2,\text{ when } \operatorname{Im}(z)<0,
\end{array} \right.
\end{aligned}
\end{equation}
where 
\begin{multicols}{2}
\begin{itemize}
    \item $A^\pm=\mp1+ib_1^\pm;$
    %\item $a_1^+=-1;$
    \item $b^\pm_1=\dfrac{p^*\mp\tau(\tau\pm2)}{2\tau}$;
    \item $\tau\in\mathbb{R}\setminus\{0\};$
    \item $p^*=4+8\displaystyle\sum_{k=1}^5e^{k\pi}+4e^{6\pi}.$
\end{itemize}
\end{multicols}
Let us prove that  there exist real parameters $s_1,s_2$ and $s_3$ small enough such that it has 3 limit cycles bifurcating from the infinity. This will be done by reducing it to the system studied in  Proposition \ref{imp_prop} by a suitable change of variables.

   Take the conformal change of variables $w={1}/{z}$. Thus, using that $\dot z=-\dot{z}/z^2,$ the piecewise system \eqref{pw_3cycles}  is transformed into the new piecewise holomorphic system
       \begin{equation}\label{pw_fam3_ex}
\begin{aligned}
\left\{\begin{array}{l}
\dot{w}=(i+1+s_1)w+(A^++s_2)w^2, \text{ when } \operatorname{Im}(w)> 0,\\[5pt]
\dot{w}=(i-1)w+(A^-+is_3)w^2,\text{ when } \operatorname{Im}(w)<0.
\end{array} \right.
\end{aligned}
\end{equation}
It behaves in a neighborhood of the origin like system \eqref{pw_3cycles} near infinity.  

From Proposition \ref{imp_prop}, we know that the origin of system \eqref{pw_fam3_ex} is a weak focus of order 4 when $s_1=s_2=s_3=0$ and there exist real parameters $s_1,s_2$ and $s_3$ small enough such that 3 limit cycles bifurcating from the origin. These 3 limit cycles are transformed, via  conformal map $z=1/w,$ into 3 limit cycles of the original system near infinity. This completes our first  proof.

\bigskip

Using the above approach it is also possible to provide numerical evidences of a concrete example of a piecewise linear holomorphic system with  3 limit cycles as was done in \cite{MR4189023}. Here, to give a proof of the existence of these 3 limit cycles for a specific system in our second proof we will use  Poincaré-Miranda theorem, as it was done by Gasull et al. in \cite{Armengol} for the same question but for piecewise linear systems.

\subsection{Second proof of Theorem \ref{teod}}\label{PM_3CL}
Let us prove that the PWHS
       \begin{equation}\label{pw_2cycles}
\begin{aligned}
\left\{\begin{array}{l}
\dot{z}=\left(i+\frac{3}{8}\right)z-\frac{14333}{2000}+i\frac{1159}{1000}, \text{ when } \operatorname{Im}(z)> 0,\\[5pt]
\dot{z}=\left(i-\frac{1}{5}\right)z-\frac{51}{50}+i\frac{1}{250},\text{ when } \operatorname{Im}(z)<0,
\end{array} \right.
\end{aligned}
\end{equation}
has at least 3 nested limit cycles surrounding the real focus $(-1/5,-49/50)$.

Denote by $z^\pm(t,x)$ the flow associated to system \eqref{pw_2cycles}, where $x=x+i0>0$ is the initial condition. The solution of the first equation (resp. second equation) of \eqref{pw_2cycles} is $z^+(t)=x^+(t)+iy^+(t)$ (resp. $z^-(t)=x^-(t)+iy^-(t)$), where
$$\begin{array}{rcl}
x^-(t)&=&-\frac{1}{5} + ((x+\frac{1}{5}) \cos(t) -\frac{49}{50} \sin(t))e^{-\frac{t}{5}}; \vspace{0.3cm}\\
y^-(t)&=&-\frac{49}{50} + ((x+\frac{1}{5})\sin(t) +\frac{49}{50} \cos(t))e^{-\frac{t}{5}};\vspace{0.3cm}\\
x^+(t)&=&\frac{67}{50} + ((x-\frac{67}{50}) \cos(t) -\frac{833}{125} \sin(t))e^{\frac{3}{8}t};\vspace{0.3cm}\\
y^+(t)&=&-\frac{833}{125} + ((x- \frac{67}{50})\sin(t) +\frac{833}{125} \cos(t))e^{\frac{3}{8}t}.\vspace{0.3cm}\\
\end{array}$$
Three equations that impose the existence of a periodic orbit are given by: $y^-(-u)=0,$ $y^+(v)=0$ and $x^-(-u)-x^+(v)=0,$ where $u>0$ and $v>0$ are the flight times and we have taken $-u$ on $\operatorname{Im}(z)<0$ since the field is reversed, and $v>0$ on $\operatorname{Im}(z)>0.$
Therefore, 
\begin{equation}\label{ee1}
    -\frac{49}{50}+\left(-\left(x+\frac{1}{5}\right)\sin(u)+\frac{49}{50}\cos(u)\right)e^{\frac{u}{5}}=0;
\end{equation}
         \begin{equation}\label{ee2}
   -\frac{833}{125}+\left(\left(x-\frac{67}{50}\right)\sin(v)+\frac{833}{125}\cos(v)\right)e^{\frac{3}{8}v}=0;        
         \end{equation}
       \begin{equation}\label{ee3}
         -\frac{77}{50}+\left(\left(x+\frac{1}{5}\right)\cos(u)+\frac{49}{50}\sin(u)\right)e^{\frac{u}{5}}-\left(\left(x-\frac{67}{50}\right)\cos(v)-\frac{833}{125}\sin(v)\right)e^{\frac{3}{8}v}=0.  
       \end{equation}
From \eqref{ee1}, we obtain that 
\begin{equation}\label{ee4}
    x=-\dfrac{10\sin(u)e^{\frac{u}{5}}-49\cos(u)e^{\frac{u}{5}}+49}{50\sin(u)e^{\frac{u}{5}}}.
\end{equation}
Substituting \eqref{ee4} in \eqref{ee2} and \eqref{ee3}, we get that $h_i(u,v)=e^{\frac{3}{8}v}\cdot e_i(u,v)$ for $i=1,2$, where  
\begin{equation}\label{ee5}
\begin{array}{rcl}
e_1(u,v)&:=&a(u)\cos(v)+b(u)\sin(v)+c(u)e^{-\frac{3}{8}v}=0,\vspace{0.3cm}\\
e_2(u,v)&:=&d(u)\cos(v)+e(u)\sin(v)+f(u)e^{-\frac{3}{8}v}=0,\vspace{0.3cm}\\
\end{array}
\end{equation}
and
$$
\begin{array}{ll}
a(u):=1666\sin(u); &
b(u):=245(\cos(u)-e^{-\frac{u}{5}})-385\sin(u);\vspace{0.3cm}\\
c(u):=-1666\sin(u);&
d(u):=-245(\cos(u)-e^{-\frac{u}{5}})+385\sin(u);\vspace{0.3cm}\\
e(u):=1666\sin(u); &
f(u):=-385\sin(u)-245\cos(u)+245e^{\frac{u}{5}}.\vspace{0.3cm}\\
\end{array}$$
 	\begin{minipage}[t]{.55\textwidth}
\raggedright
\begin{figure}[H]
	\begin{center}
		\begin{overpic}[scale=0.43]{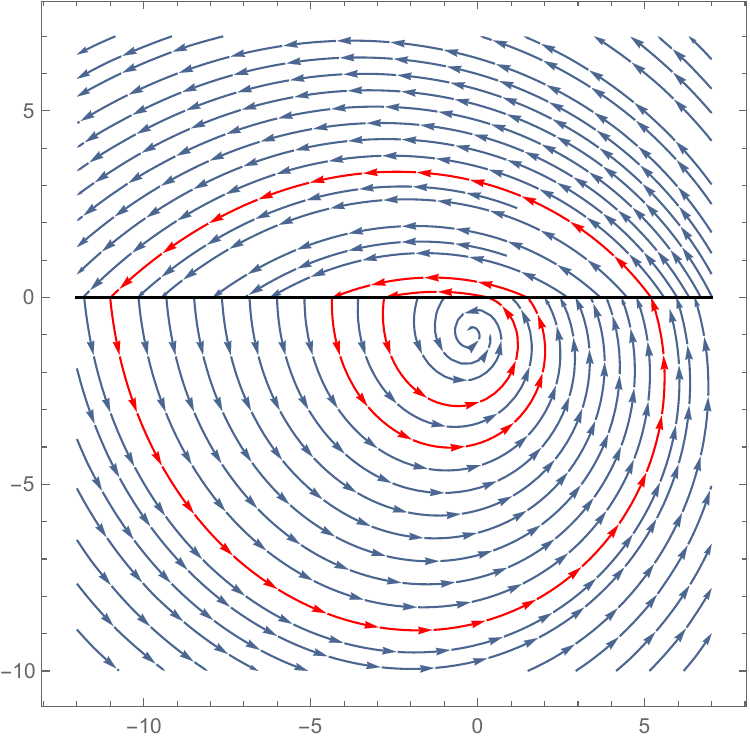}
		%\begin{overpic}[grid,tics=5,width=5cm]{limit_cycle_1.pdf}		
        \put(-7,31){$\Sigma^-$}
        \put(-7,75){$\Sigma^+$}
		\put(102,58){$\Sigma$}
		\end{overpic}
		\caption{\footnotesize{Phase portrait of PWHS \eqref{pw_2cycles}. The red trajectories are the 3 nested limit cycles of \eqref{pw_2cycles} surrounding the focus $(-{1}/{5}, -{49}/{50})$.}}
	\label{limit_cycle_2}
	\end{center}
	\end{figure}
\end{minipage}
\begin{minipage}[t]{.55\textwidth}
\raggedright
\begin{figure}[H]
	\begin{center}
		\begin{overpic}[scale=0.42]{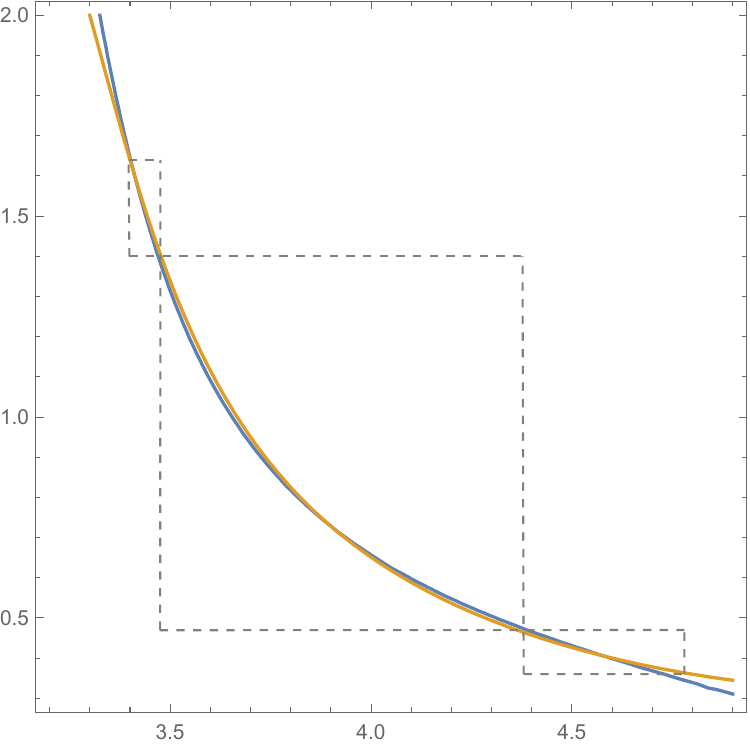}
		%\begin{overpic}[grid,tics=5,width=5cm]{limit_cycle_1.pdf}		
        \put(50,50){$B_2$}
        \put(7,75){$B_1$}
		\put(89,18){$B_3$}
		\end{overpic}
		\caption{\footnotesize{Boxes $B_1$, $B_2$, $B_3$. The blue and orange curves correspond to $e_1(u,v)$ and $e_2(u,v),$ respectively.}}
	\label{fig_boxes}
	\end{center}
	\end{figure}
\end{minipage}

\bigskip

By solving system \eqref{ee5} numerically, we obtain the 3 approximate solutions: $$(3.411939, 1.599505),\quad (3.896572,0.731139),\quad \text{and}\quad (4.583524,0.404161).$$ In what follows, we prove that near them there exist actual solutions by using the Poincaré--Miranda theorem.

Now, consider the following 3 boxes, each one of them containing one of the approximate solutions:
$$B_1=\left[\frac{17}{5},\frac{139}{40}\right]\times\left[\frac{7}{5},\frac{41}{25}\right];\quad B_2=\left[\frac{139}{40},\frac{219}{50}\right]\times\left[\frac{47}{100},\frac{7}{5}\right];\quad B_3=\left[\frac{219}{50},\frac{239}{50}\right]\times\left[\frac{9}{25},\frac{47}{100}\right].$$
They will be used to apply the Poincaré--Miranda theorem for the functions $e_1$ and $e_2,$ see Figure \ref{fig_boxes}.

We give some details in box $B_1$, in the other boxes it is done in an analogous way. We claim that $e_1(17/5,v)<0,$ for all $v\in[7/5,41/25].$ Indeed, 
$$e_1\left(\frac{17}{5},v\right)=a\left(\frac{17}{5}\right)\cos(v)+b\left(\frac{17}{5}\right)\sin(v)+c\left(\frac{17}{5}\right)e^{-\frac{v}{5}},$$
which is the form of the function $h$ given in Lemma \ref{lemma2}, where
$$
\begin{array}{ll}
A=a\left(\frac{17}{5},v\right)=1666\sin\left(\frac{17}{5}\right); &
B=b\left(\frac{17}{5},v\right)=245(\cos\left(\frac{17}{5}\right)-e^{-\frac{17}{25}})-385\sin\left(\frac{17}{5}\right); \vspace{0.3cm}\\
C=0;&
D=c\left(\frac{17}{5},v\right)=-1666\sin\left(\frac{17}{5}\right).
\end{array}$$
Applying Lemma \ref{lemma2} for $n=7$, we obtain its Taylor series expansion around $u=0$ and $\overline{m}\approx 0.0170742$. Then, maximizing said polynomial by a polynomial $p^+_{n,k}$ with rational coefficients of degree $7$ given by 
$$\begin{array}{rcl}
p^+_{n,k}(v)&=&\displaystyle\sum_{j=0}^n a_j^+v^j+Mv^{n+1},\vspace{0.3cm}\\
     &=& \frac{1}{100000}-\frac{10556317}{25000}v+\frac{24279999}{100000}v^2+\frac{2001273}{50000}v^3-\frac{4347}{250}v^4-\frac{110733}{50000}v^5\vspace{0.3cm}\\
     &&+\frac{29647}{50000}v^6+\frac{2601}{50000}v^7+Mv^8,
\end{array}$$ 
where we have used that
$a_j^+=trunc(a_j\cdot10^{k})\cdot10^{-k}+10^{-k}\in\mathbb{Q},$ $k=5$
and that $M=1/58$ is an upper bound of the right-hand side expression at \eqref{M}, i.e. $\overline{m}< M$.

From the Sturm's sequences of $p^+_{n,k}$, we obtain that it has no roots in $[7/5,41/25]$ and $e_1(17/5,v)<p^+_{n,k}<0,$ for all $v\in[7/5,41/25].$ 

To prove that $e_1(139/40,v)>0$ for all $v\in[7/5,41/25],$ we take $n=6$ and $k=2$ and we do the same previous analysis but we minorize $e_1$ by a polynomial with rational coefficients given by 
$$\begin{array}{rcl}
p^-_{n,k}(v)&=&\displaystyle\sum_{j=0}^n a_j^-v^j-Mv^{n+1},\vspace{0.3cm}\\
     &=&-\frac{1}{100}-\frac{1729}{4}v+\frac{31093}{100}v^2+\frac{829}{5}v^3-\frac{2227}{100}v^4-\frac{97}{50}v^5+\frac{37}{50}v^6-Mv^7,

\end{array}
$$
where $a_j^-=trunc(a_j\cdot10^{k})\cdot10^{-k}-10^{-k}\in\mathbb{Q}$ and $M=1/6.$ Using the Sturm's sequences of $p^-_{n,k}$, we conclude that it has no roots in $[7/5,41/25]$ and $0<p^-_{n,k}<e_1(139/40,v),$ for all $v\in[7/5,41/25].$

In the following tables we summarize the results obtained using the previous analysis in each of the boxes $B_i$, for $i=1,2,3$.

\bigskip

    \begin{equation*}\label{table0}
\begin{array}{||c|| c |c| c | c|c |c||}
\hline
B_1&\textit{Face}&	u=17/5 &u=139/40 & v=7/5 &v=41/25 \\
\hline\hline
&\textit{Function sign}&e^{\frac{3}{8}v}\cdot e_1<0&e^{\frac{3}{8}v}\cdot e_1>0 &e_2>0&e_2<0 \\
\hline
&\textit{Polynomial}&p^+_{n,k}	&p^-_{n,k} &p^-_{n,k}&p^+_{n,k}\\
\hline
&\textit{Parameters}:(n,k)&(7,5)	&(6,2) &(15,13)&(15,13) \\
\hline
&M&1/58	&1/6 &35/26\cdot10^{-10}&11/8\cdot10^{-10} \\
\hline
\end{array}
\end{equation*}
\bigskip

    \begin{equation*}\label{table1}
\begin{array}{||c|| c |c| c | c|c |c||}
\hline
B_2&\textit{Face}&	u=139/40 &u=219/50 & v=47/100 &v=7/5 \\
\hline\hline
&\textit{Function sign}&e_2>0&e_2<0 &e^{\frac{3}{8}v}\cdot e_1<0&e^{\frac{3}{8}v}\cdot e_1>0 \\
\hline
&\textit{Polynomial}&p^-_{n,k}	&p^+_{n,k} &p^+_{n,k}&p^-_{n,k}\\
\hline
&\textit{Parameters}:(n,k)&(8,5)	&(8,5) &(14,11)&(15,11) \\
\hline
&M&1/469	&1/206 &9/5\cdot10^{-10}&15/14\cdot10^{-10} \\
\hline
\end{array}
\end{equation*}
\bigskip

    \begin{equation*}\label{table2}
\begin{array}{||c|| c |c| c | c|c |c||}
\hline
B_3&\textit{Face}&	u=219/50 &u=239/50 & v=9/25 &v=47/100 \\
\hline\hline
&\textit{Function sign}&e^{\frac{3}{8}v}\cdot e_1<0&e^{\frac{3}{8}v}\cdot e_1>0 &e_2>0&e_2<0 \\
\hline
&\textit{Polynomial}&p^+_{n,k}	&p^-_{n,k} &p^-_{n,k}&p^+_{n,k}\\
\hline
&\textit{Parameters}:(n,k)&(6,3)	&(6,3) &(14,12)&(14,11) \\
\hline
&M&7/20	&2/5 &37/4\cdot10^{-10}&11/10\cdot10^{-10} \\
\hline
\end{array}
\end{equation*}

\bigskip

The information given in the table corresponding to $B_1$ together with the Poincaré--Miranda theorem \ref{PM} allows us to conclude the existence of a positive zero of the system \eqref{ee5} in this box. Doing the same analysis in $B_2$ and $B_3$ allows us to conclude the existence of a positive zero of the system \eqref{ee5} in each box. From these 3 solutions, 3 nested limit cycles arise surrounding the real focus $(-1/5,-49/50)$, see Figure \ref{limit_cycle_2}.

The proof of the hyperbolicity of the limit cycles follows the same ideas as \cite{Armengol} and \cite{320978dfa23348229b3918bd894574f5}.

 %%%%%%%%%%%%%%%%%%%%%%%%%%%5
\section*{Appendix: Proof of Proposition \ref{main:prop}}\label{appendix}
    In what follows, we will use the following formulas. Let $\alpha$ and $\beta$ be complex numbers, then
$$\begin{array}{ll}
2\operatorname{Re}(\alpha)\operatorname{Re}(\beta)=\operatorname{Re}[\alpha\beta+\overline{\alpha}\beta];&\operatorname{Re}(\alpha)=\operatorname{Im}(i\alpha);\\
2\operatorname{Im}(\alpha)\operatorname{Im}(\beta)=\operatorname{Re}[-\alpha\beta+\overline{\alpha}\beta];&\operatorname{Re}(i\alpha)=-\operatorname{Im}(\alpha).\\
2\operatorname{Re}(\alpha)\operatorname{Im}(\beta)=\operatorname{Im}[\alpha\beta+\overline{\alpha}\beta];\\
\end{array}$$ 
         
         First, we compute $\omega_1(\pi)$. Recall that $F_1(z)=(\lambda+i)z,$ thus $S_1(\theta)=\lambda+i.$ Consider the differential equation \eqref{smooth_eq1} in polar coordinates, which is given as \eqref{smooth_eq_3}. Direct substitution shows that $\omega_1(\theta)$ satisfies
         $$\omega'_1(\theta)=\frac{\operatorname{Re}(S_1(\theta))}{\operatorname{Im}(S_1(\theta))}\omega_1(\theta)= \lambda \omega_1(\theta),\quad \omega_1(0)=0.$$
     Hence
     $\omega_1(\theta)=e^{\lambda\theta}-1
     $ and  $\omega_1(\pi)=e^{\lambda \pi}-1.$
     Now, we compute $\omega_2(\pi)$ and $\omega_3(\pi).$  Recall that equation \eqref{smooth_eq_3} in polar coordinates $R,\theta$ writes as
     \begin{equation}\label{smooth_eq_5}
    \dfrac{dR}{d\theta}=\dfrac{\lambda R+\sum_{k=2}^\infty R^k\operatorname{Re}(S_k(\theta))}{1+\sum_{k=2}^\infty R^{k-1}\operatorname{Im}(S_k(\theta))}=\sum_{k=1}^\infty T_k(\theta)R^k,
\end{equation}
where
$$\begin{array}{rcl}
    T_1 &=& \lambda; \\
    T_2 &=&  \operatorname{Re}(S_2)-\lambda \operatorname{Im}(S_2); \\
    T_3 &=& \operatorname{Re}(S_3)-\operatorname{Re}(S_2)\operatorname{Im}(S_2)-\lambda\operatorname{Im}(S_3)+\lambda(\operatorname{Im}(S_2))^2;\\
    T_4&=& \operatorname{Re}(S_2)(\operatorname{Im}(S_2)^2-\operatorname{Im}(S_3))-\operatorname{Im}(S_2) \operatorname{Re}(S_3)+\operatorname{Re}(S_4)\\
    & &+\lambda(-\operatorname{Im}(S_2)^3 + 2 \operatorname{Im}(S_2)\operatorname{Im}(S_3)-\operatorname{Im}(S_4));\\
    T_5&=&\operatorname{Re}(S_2)(-\operatorname{Im}(S_2)^3+2\operatorname{Im}(S_2)\operatorname{Im}(S_3)-\operatorname{Im}(S_4))\\
    & &+\operatorname{Re}(S_3)(\operatorname{Im}(S_2)^2-\operatorname{Im}(S_3))-\operatorname{Im}(S_2)\operatorname{Re}(S_4)+\operatorname{Re}(S_5)\\
    & &+\lambda(\operatorname{Im}(S_2)^4-3\operatorname{Im}(S_2)^2\operatorname{Im}(S_3)+\operatorname{Im}(S_3)^2+2\operatorname{Im}(S_2)\operatorname{Im}(S_4)-\operatorname{Im}(S_5)).
\end{array}$$
Consider the change of variables $r=Re^{-\lambda\theta}.$ Thus equation \eqref{smooth_eq_5} is converted into 
\begin{equation}\label{smooth_eq_7}
    \dfrac{dr}{d\theta}=\sum_{k=2}^\infty R_k(\theta)r^k,
\end{equation}
where $R_k(\theta)=T_k(\theta)e^{(k-1)\lambda\theta},$ for each $k\geq 2.$ Therefore, 
$$\begin{array}{rcl}
    R_2 &=& e^{\lambda\theta}[\operatorname{Re}(S_2)-\lambda \operatorname{Im}(S_2)]; \\
    R_3 &=& e^{2\lambda\theta}[\operatorname{Re}(S_3)-\operatorname{Re}(S_2)\operatorname{Im}(S_2)-\lambda\operatorname{Im}(S_3)+\lambda(\operatorname{Im}(S_2))^2];\\
     R_4&=& e^{3\lambda\theta}[\operatorname{Re}(S_2)(\operatorname{Im}(S_2)^2-\operatorname{Im}(S_3))-\operatorname{Im}(S_2) \operatorname{Re}(S_3)+\operatorname{Re}(S_4)\\
    & &+\lambda(-\operatorname{Im}(S_2)^3 + 2 \operatorname{Im}(S_2)\operatorname{Im}(S_3)-\operatorname{Im}(S_4))];\\
    R_5&=&e^{4\lambda\theta}[\operatorname{Re}(S_2)(-\operatorname{Im}(S_2)^3+2\operatorname{Im}(S_2)\operatorname{Im}(S_3)-\operatorname{Im}(S_4))\\
    & &+\operatorname{Re}(S_3)(\operatorname{Im}(S_2)^2-\operatorname{Im}(S_3))-\operatorname{Im}(S_2)\operatorname{Re}(S_4)+\operatorname{Re}(S_5)\\
    & &+\lambda(\operatorname{Im}(S_2)^4-3\operatorname{Im}(S_2)^2\operatorname{Im}(S_3)+\operatorname{Im}(S_3)^2+2\operatorname{Im}(S_2)\operatorname{Im}(S_4)-\operatorname{Im}(S_5))].
\end{array}$$
%Consequently, from  Proposition, when we consider its solution in the form \eqref{smooth_eq_3}
Let $r$ be the solution of \eqref{smooth_eq_7} such that $r(0,\rho)=\rho.$ Then
\begin{equation*}%\label{polar_sol_1}
    r(\theta,\rho)-\rho=\sum_{k=2}^\infty u_k(\theta)\rho^k,\,\text{where}\,u_k(0)=0,\,\text{for}\, k\geq 2.
\end{equation*}
Since $r=Re^{-\lambda\theta},$ then $$
\begin{array}{rcl}
R(\theta,s)&=&r(\theta,s)e^{\lambda\theta},\\
&=&\left[s+\sum_{k=2}^\infty u_k(\theta)s^k\right]e^{\lambda\theta}.\\
\end{array}
$$
Thus, 
\begin{equation}\label{wurel}
    R(\theta,s)=s+\sum_{k=1}^\infty \omega_k(\theta)s^k=\left[s+\sum_{k=2}^\infty u_k(\theta)s^k\right]e^{\lambda\theta}.
\end{equation}
By equation \eqref{wurel}, we get that $\omega_k(\pi)=e^{\lambda\pi}u_k(\pi),$ for all $k\geq 2.$
In what follows, we compute the $u_i$, for each $i=2,\cdots,5$. We start calculating $u_2$. From Proposition \ref{coef_R}, $u_2(\theta)=\widetilde{R}_2.$ Thus, 
$$
\begin{array}{rcl}
\widetilde{R_2}(\theta)
     &=& \displaystyle\operatorname{Re}\int_{0}^\theta S_2(\varphi)e^{\lambda\varphi} d\varphi-\lambda\operatorname{Im}\int_{0}^\theta S_2(\varphi)e^{\lambda\varphi} d\varphi\\
     &=& \displaystyle\operatorname{Re}\left[\int_{0}^\theta(1+\lambda i) S_2(\varphi)e^{\lambda\varphi} d\varphi\right]\\
      &=& \displaystyle\operatorname{Re}\underbrace{\left[\frac{(1+\lambda i) A(e^{(\lambda+i)\theta}-1)}{i+\lambda}\right]}_{\substack{\kappa}}.
\end{array}
$$
Hence, $$\omega_2(\pi)=e^{\lambda\pi}u_2(\pi)=e^{\lambda\pi}(-e^{\lambda\pi}-1)\displaystyle\operatorname{Re}\left[\frac{(1+\lambda i) A}{i+\lambda}\right].$$
Now, we find $u_3$. From Proposition \ref{coef_R},  $u_3(\theta)=(\widetilde{R}_2)^2+\widetilde{R}_3.$ Thus, 
$$\begin{array}{rcl}
 \widetilde{R_3}(\theta)
     &=& \operatorname{Re}\displaystyle\int_{0}^{\theta} e^{2\lambda\varphi}\left[S_3(\varphi)+\frac{\lambda}{2}(-S_2^2(\varphi)+S_2(\varphi)\overline{S_2}(\varphi))\right]d\varphi\vspace{0.3cm}\\
     & &-\operatorname{Im}\displaystyle\int_{0}^\theta e^{2\lambda\varphi}\left[\lambda S_3(\varphi)+\frac{1}{2}S_2^2(\varphi)\right]d\varphi\vspace{0.3cm}\\
&=&\operatorname{Re}\displaystyle\int_{0}^{\theta}e^{2\lambda\varphi}(1+\lambda i)S_3(\varphi)d\varphi+\dfrac{\lambda}{2}\operatorname{Re}\displaystyle\int_{0}^{\theta}e^{2\lambda\varphi}S_2(\varphi)\overline{S_2}(\varphi)d\varphi\vspace{0.3cm}\\ 
&&-\dfrac{1}{2}\operatorname{Im}\displaystyle\int_{0}^{\theta}e^{2\lambda\varphi}(1+\lambda i)(S_2(\varphi))^2 d\varphi\vspace{0.3cm}\\
     &=&\operatorname{Re}\underbrace{\left[\dfrac{(1+\lambda i)B(e^{2\theta(i+\lambda)}-1)}{2(i+\lambda)}+\dfrac{A\overline{A}(e^{2\lambda\theta}-1)}{4}\right]}_{\substack{\alpha}}\vspace{0.3cm}\\
     &&-\operatorname{Im}\underbrace{\left[ \dfrac{(1+\lambda i)A^2(e^{2\theta(i+\lambda)}-1)}{4(i+\lambda)}\right]}_{\substack{\beta}}.\\
    % &=&\operatorname{Re}(\alpha)-\operatorname{Im}(\beta)
  \end{array}
$$
Therefore, 
$$\begin{array}{rcl}
\omega_3(\pi)&=&e^{\lambda\pi}u_3(\pi)\\
&=& e^{-\lambda\pi}(\omega_2(\pi))^2+e^{\lambda\pi}(e^{2\lambda\pi}-1)\left\{\operatorname{Re}\left[\dfrac{(1+\lambda i)B}{2(i+\lambda)}+\dfrac{ A\overline{A}}{4}\right]-\operatorname{Im}\left[\dfrac{(1+\lambda i)A^2}{4(i+\lambda)}\right]\right\}.
\end{array}$$
Now, we compute $u_4$. From Proposition \ref{coef_R},  $u_4(\theta)=(\widetilde{R}_2)^3+3\widetilde{R}_2\widetilde{R}_3-\widetilde{\widetilde{R}_3R_2}+\widetilde{R}_4.$ Thus, 

$$\begin{array}{rcl}
\widetilde{\widetilde{R_3}R_2}(\pi)&=&\dfrac{1}{2}\operatorname{Re}\displaystyle\int_{0}^{\pi}e^{\lambda\varphi}(\alpha+\overline{\alpha})S_2(\varphi)d\varphi+\dfrac{\lambda}{2}\operatorname{Re}\displaystyle\int_{0}^{\pi}e^{\lambda\varphi}(-S_2(\varphi)+\overline{S_2}(\varphi))\beta d\varphi\\ 
&&-\dfrac{\lambda}{2}\operatorname{Im}\displaystyle\int_{0}^{\pi}e^{\lambda\varphi}(\alpha+\overline{\alpha})S_2(\varphi)d\varphi-\dfrac{1}{2}\operatorname{Im}\displaystyle\int_{0}^{\pi}e^{\lambda\varphi}(S_2(\varphi)+\overline{S_2}(\varphi))\beta d\varphi\\
 &=&\dfrac{1}{2}(-e^{3\pi\lambda}-1)\left\{\operatorname{Re}\left[\eta_1-\lambda\eta_5\right]-\operatorname{Im}\left[\lambda\eta_1+\eta_3\right]\right\}\vspace{0.3cm}\\
 &&-\dfrac{1}{2}(-e^{\pi\lambda}-1)\left\{\operatorname{Re}\left[\eta_2-\lambda\eta_6\right]-\operatorname{Im}\left[\lambda\eta_2+\eta_4\right]\right\},\\
%&=&\frac{1}{2}\eta_1(\pi)-\frac{\lambda}{2}\eta_2(\pi)-\frac{1}{2}\eta_3(\pi)+\frac{\lambda}{2}\eta_4(\pi)\\
\end{array}$$
where
$$\begin{array}{ll}
\eta_1=\dfrac{(1+\lambda i)AB}{6(i+\lambda)^2}+\dfrac{(1-\lambda i)A\overline{B}}{2(-i+\lambda)(3\lambda-i)}+\dfrac{A^2\overline{A}}{2(i+3\lambda)};& \eta_4=\dfrac{(1+\lambda i)A^3}{4(i+\lambda)^2}+\dfrac{(1+\lambda i)A^2\overline{A}}{4(1+\lambda^2)};\vspace{0.3cm}\\
\eta_2=\dfrac{(1+\lambda i)AB}{2(i+\lambda)^2}+\dfrac{(1-\lambda i)A\overline{B}}{2(1+\lambda^2)}+\dfrac{A^2\overline{A}}{2(i+\lambda)};& \eta_5=\dfrac{(1+\lambda i)A^3}{12(i+\lambda)^2}-\dfrac{(1+\lambda i)A^2\overline{A}}{4(i+\lambda)(3\lambda+i)};\vspace{0.3cm}\\
\eta_3=\dfrac{(1+\lambda i)A^3}{12(i+\lambda)^2}+\dfrac{(1+\lambda i)A^2\overline{A}}{4(i+\lambda)(3\lambda+i)}; &\eta_6=\dfrac{(1+\lambda i)A^3}{4(i+\lambda)^2}-\dfrac{(1+\lambda i)A^2\overline{A}}{4(1+\lambda^2)}.\vspace{0.3cm}\\
\end{array}$$

$$\begin{array}{rcl}
\widetilde{R_4}(\theta)&=&
 %\displaystyle\int_{0}^{\psi}R_4(\psi) d\psi&=&
 %\displaystyle\int_{0}^{\pi}e^{3\lambda\varphi}[\operatorname{Re}%(S_2(\varphi))(\operatorname{Im}(S_2(\varphi))^2-\operatorname{Im}(S_3)(\varphi))\vspace{0.3cm}\\
 %&&-\operatorname{Im}(S_2(\varphi)) \operatorname{Re}(S_3(\varphi))+\operatorname{Re}(S_4(\varphi))\vspace{0.3cm}\\
 %&&+\lambda(-\operatorname{Im}(S_2)^3 + 2 \operatorname{Im}(S_2)\operatorname{Im}(S_3)-\operatorname{Im}(S_4))]d\varphi\vspace{0.3cm}\\ %%%%%%%%%%%%%%%%%%%%%%%
 \dfrac{1}{4}\operatorname{Re}\displaystyle\int_{0}^{\theta}e^{3\lambda\varphi}\left[-(S_2(\varphi))^3+(\overline{S_2}(\varphi))^2S_2(\varphi)+4\lambda(-S_2(\varphi)S_3(\varphi)\right.\vspace{0.3cm}\\
&&+\left.\overline{S_2}(\varphi)S_3(\varphi))+4S_4(\varphi)\right]d\varphi-\dfrac{1}{2}\operatorname{Im}\displaystyle\int_{0}^{\theta}e^{3\lambda\varphi}\left[2S_2(\varphi)S_3(\varphi)\right.\vspace{0.3cm}\\
 &&+\overline{S_2}(\varphi)S_3(\varphi)+\overline{S_3}(\varphi)S_2(\varphi)+\frac{\lambda}{2}(-(S_2(\varphi))^3+2\overline{S_2}(\varphi)(S_2(\varphi))^2\vspace{0.3cm}\\
 & &\left.-(\overline{S_2}(\varphi))^2S_2(\varphi)+4S_4(\varphi))\right]d\varphi\vspace{0.3cm}\\
 &=&-\dfrac{1}{4}\operatorname{Re}\displaystyle\int_{0}^{\theta}e^{3\lambda\varphi}(1+\lambda i)(S_2(\varphi))^3 d\varphi+\operatorname{Re}\displaystyle\int_{0}^{\theta}e^{3\lambda\varphi}(1+\lambda i)S_4(\varphi) d\varphi\\ &&+\operatorname{Re}\displaystyle\int_{0}^{\theta}e^{3\lambda\varphi}(i-\lambda)S_2(\varphi)S_3(\varphi) d\varphi+\dfrac{1}{2}\operatorname{Re}\displaystyle\int_{0}^{\theta}e^{3\lambda\varphi}(2\lambda-i)\overline{S_2}(\varphi)S_3(\varphi) d\varphi\\
 &&+\dfrac{1}{4}\operatorname{Im}\displaystyle\int_{0}^{\theta}e^{3\lambda\varphi}(i+\lambda)(\overline{S_2}(\varphi))^2S_2(\varphi) d\varphi-\dfrac{1}{2}\operatorname{Im}\displaystyle\int_{0}^{\theta}e^{3\lambda\varphi}\overline{S_3}(\varphi)S_2(\varphi) d\varphi\\
  &&-\dfrac{\lambda}{2}\operatorname{Im}\displaystyle\int_{0}^{\theta}e^{3\lambda\varphi}\overline{S_2}(\varphi)(S_2(\varphi))^2 d\varphi\\
%        \end{array}$$
%     $$\begin{array}{rcl}
  &=&\operatorname{Re}\left[ \dfrac{(1+i\lambda)(e^{3\theta(i+\lambda)}-1)}{3(i+\lambda)}\left(-\dfrac{A^3}{4}+C+iAB\right)+\dfrac{(2\lambda-i)(e^{\theta(3\lambda+i)}-1)\overline{A}B}{2(3\lambda+i)}\right]\vspace{0.3cm}\\
     &&+\operatorname{Im}\left[\dfrac{(e^{\theta(3\lambda-i)}-1)}{2(3\lambda-i)}\left(\dfrac{(i+\lambda)(\overline{A})^2A}{2}-A\overline{B}\right)- \dfrac{(e^{\theta(3\lambda+i)}-1)\lambda\overline{A}A^2}{2(3\lambda+i)}\right]\vspace{0.3cm}\\
     &=&\operatorname{Re}[\underbrace{(e^{3\theta(i+\lambda)}-1)\gamma_1+(e^{\theta(3\lambda+i)}-1)\gamma_2}_{\substack{\gamma_R} }]+\operatorname{Im}[\underbrace{(e^{\theta(3\lambda-i)}-1)\gamma_3+(e^{\theta(3\lambda+i)}-1)\gamma_4}_{\substack{\gamma_I} }],
\end{array}$$
where
$$
\begin{array}{ll}
     \gamma_1= \dfrac{(1+i\lambda)}{3(i+\lambda)}\left(-\dfrac{A^3}{4}+C+iAB\right);&\gamma_3= \dfrac{1}{2(3\lambda-i)}\left(\dfrac{(i+\lambda)(\overline{A})^2A}{2}-A\overline{B}\right);\\
      \gamma_2=\dfrac{(2\lambda-i)\overline{A}B}{2(3\lambda+i)};&   \gamma_4= - \dfrac{\lambda\overline{A}A^2}{2(3\lambda+i)}.
\end{array}
$$
Consequently,
$$\begin{array}{rcl}
\omega_4(\pi)&=&e^{\lambda\pi}u_4(\pi)\\
&=&e^{\lambda\pi}\left(\left(u_2(\pi)\right)^3+3u_2(\pi)(u_3(\pi)-(u_2(\pi))^2)-\widetilde{\widetilde{R_3}R_2}(\pi)+\widetilde{R_4}(\pi)\right)\\
&=&-2e^{-2\lambda\pi}\left(\omega_2(\pi)\right)^3+3e^{-\lambda\pi}\omega_2(\pi)\omega_3(\pi)\\
&&-\dfrac{1}{2}e^{\lambda\pi}(-e^{3\pi\lambda}-1)\left\{\operatorname{Re}\left[\eta_1-\lambda\eta_5-2(\gamma_1+\gamma_2)\right]-\operatorname{Im}\left[\lambda\eta_1+\eta_3+2(\gamma_3+\gamma_4)\right]\right\}\vspace{0.3cm}\\
 &&+\dfrac{1}{2}e^{\lambda\pi}(-e^{\pi\lambda}-1)\left\{\operatorname{Re}\left[\eta_2-\lambda\eta_6\right]-\operatorname{Im}\left[\lambda\eta_2+\eta_4\right]\right\}.\\
\end{array}$$

Finally, we compute $u_5$. From Proposition \ref{coef_R},  $u_5(\theta)=(\widetilde{R}_2)^4+5(\widetilde{R}_2)^2\widetilde{R}_3+\widetilde{(\widetilde{R}_2)^2R_3}-2\widetilde{R}_2\widetilde{\widetilde{R}_3R_2}+\frac{3}{2}(\widetilde{R}_3)^2+4\widetilde{R}_2\widetilde{R}_4-2\widetilde{\widetilde{R}_4R_2}+\widetilde{R}_5.$ Thus, 

$$\begin{array}{rcl}
\widetilde{R_5}(\pi)&=&
-\dfrac{1}{8}\operatorname{Re}\displaystyle\int_{0}^{\pi}e^{4\lambda\varphi}(S_2(\varphi)+\overline{S_2}(\varphi))^2(-(S_2(\varphi))^2+\overline{S_2}(\varphi)S_2(\varphi)) d\varphi\\ &&+\dfrac{1}{2}\operatorname{Re}\displaystyle\int_{0}^{\pi}e^{4\lambda\varphi}(S_2(\varphi)+\overline{S_2}(\varphi))(-S_2(\varphi)S_3(\varphi)+\overline{S_2}(\varphi)S_3(\varphi)) d\varphi\\
&& -\dfrac{1}{2}\operatorname{Im}\displaystyle\int_{0}^{\pi}e^{4\lambda\varphi}(S_2(\varphi)+\overline{S_2(\varphi)})S_4(\varphi) d\varphi\\
 &&+\dfrac{1}{4}\operatorname{Re}\displaystyle\int_{0}^{\pi}e^{4\lambda\varphi}(S_3(\varphi)+\overline{S_3}(\varphi))(-(S_2(\varphi))^2+\overline{S_2}(\varphi)S_2(\varphi)) d\varphi\\
  && -\dfrac{1}{2}\operatorname{Im}\displaystyle\int_{0}^{\pi}e^{4\lambda\varphi}((S_3(\varphi))^2+\overline{S_3}(\varphi)S_3(\varphi))) d\varphi\\
  &&-\dfrac{1}{2}\operatorname{Im}\displaystyle\int_{0}^{\pi}e^{4\lambda\varphi}((S_4(\varphi)+\overline{S_4}(\varphi))S_2(\varphi)) d\varphi\\
  &&+\dfrac{\lambda}{8}\operatorname{Re}\displaystyle\int_{0}^{\pi}e^{4\lambda\varphi}\left[(-(S_2(\varphi))^2+\overline{S_2}(\varphi)S_2(\varphi))^2+|-(S_2(\varphi))^2+\overline{S_2}(\varphi)S_2(\varphi)|^2\right] d\varphi\\
  &&-\dfrac{3\lambda}{4}\operatorname{Im}\displaystyle\int_{0}^{\pi}e^{4\lambda\varphi}\left[-(S_2(\varphi)-\overline{S_2}(\varphi))^2S_3(\varphi)\right] d\varphi\\
  &&+\dfrac{\lambda}{2}\operatorname{Re}\displaystyle\int_{0}^{\pi}e^{4\lambda\varphi}(-(S_3(\varphi))^2+\overline{S_3}(\varphi)S_3(\varphi)) d\varphi\\
&&+\lambda\operatorname{Re}\displaystyle\int_{0}^{\pi}e^{4\lambda\varphi}(-S_2(\varphi)S_4(\varphi)+\overline{S_2}(\varphi)S_4(\varphi)) d\varphi+\operatorname{Re}\displaystyle\int_{0}^{\pi}e^{4\lambda\varphi}(1+\lambda i)S_5(\varphi) d\varphi\\ %%%%%%%%%%%%%%%%%%%%
    &=&(e^{4\pi\lambda}-1)\left\{\operatorname{Re}\left[\xi_1+\xi_3+\xi_5+\xi_8+\xi_9+\xi_{10}+\xi_{11}+\xi_{12}\right]+\operatorname{Im}\left[\xi_6+\xi_7\right]\right\}\vspace{0.3cm}\\
     &&+(-e^{4\pi\lambda}-1)\left\{\operatorname{Re}\left[\xi_2\right]-\operatorname{Im}\left[\xi_4\right]\right\},\vspace{0.3cm}\\
\end{array}$$
where
$$\begin{array}{ll}
\xi_1=\dfrac{A^4}{32(i+\lambda)}+\dfrac{A^3\overline{A}}{16(i+2\lambda)}-\dfrac{A\overline{A}^3}{16(-i+2\lambda)};&\xi_7=-\dfrac{A(2\overline{C}(i+\lambda)+C(-i+2\lambda))}{2(4+4i\lambda+8\lambda)^2};\vspace{0.3cm}\\
\xi_2=-\dfrac{A^2\overline{A}}{8(i+4\lambda)};&\xi_8=\dfrac{\lambda A}{32}\left(\dfrac{2\overline{A}^3}{i-2\lambda}+\dfrac{3A\overline{A}^2}{\lambda}+\dfrac{A^3}{i+\lambda}-\dfrac{6A^2\overline{A}}{i+2\lambda}\right);\vspace{0.3cm}\\
\xi_3=-\dfrac{(A\lambda+\overline{A}(i+\lambda))AB}{8(-1+2i\lambda)\lambda(i+\lambda)};&\xi_9=\dfrac{3B(-\overline{A}^2+i(A-3\overline{A})(A-\overline{A})\lambda+2(A-\overline{A})^2\lambda^2)}{16(i+\lambda)(i+2\lambda)};\vspace{0.3cm}\\
\xi_4=-\dfrac{D}{2}\left(-\dfrac{\overline{A}}{3i+4\lambda}-\dfrac{A}{5i+4\lambda}\right);&\xi_{10}=-\dfrac{B(B\lambda-\overline{B}(i+\lambda))}{8(i+\lambda)};\vspace{0.3cm}\\
\xi_5=\dfrac{A}{16}\left(-\dfrac{A\overline{B}}{\lambda}-\dfrac{AB}{i+\lambda}+\dfrac{2\overline{A}\overline{B}}{-i+2\lambda}+\dfrac{\overline{A}B}{i+2\lambda}\right);&\xi_{11}=\dfrac{\lambda C}{4}\left(-\dfrac{A}{i+\lambda}+\dfrac{2\overline{A}}{i+2\lambda}\right);\vspace{0.3cm}\\
\xi_6=-\dfrac{B(B\lambda+\overline{B}(i+\lambda))}{8\lambda(i+\lambda)};&\xi_{12}=\dfrac{D(1+i\lambda)}{i+\lambda}.\vspace{0.3cm}\\
\end{array}$$

$$\begin{array}{rcl}
\widetilde{\widetilde{R_4}R_2}(\pi)&=& \dfrac{1}{2}\displaystyle\int_{0}^{\varphi}e^{\lambda\theta}\operatorname{Re}[\gamma_RS_2(\varphi)+\overline{\gamma_R}S_2(\varphi)]d\varphi-\dfrac{\lambda}{2}\displaystyle\int_{0}^{\varphi}e^{\lambda\theta}\operatorname{Im}[\gamma_RS_2(\varphi)+\overline{\gamma_R}S_2(\varphi)]d\varphi\vspace{0.3cm}\\
&&+ \dfrac{1}{2}\displaystyle\int_{0}^{\varphi}e^{\lambda\theta}\operatorname{Im}[\gamma_IS_2(\varphi)+\overline{\gamma_I}S_2(\varphi)]d\varphi-\dfrac{\lambda}{2}\displaystyle\int_{0}^{\varphi}e^{\lambda\theta}\operatorname{Re}[\gamma_IS_2(\varphi)+\overline{\gamma_I}S_2(\varphi)]d\varphi\vspace{0.3cm}\\
&=&(e^{4\pi\lambda}-1)\left\{\operatorname{Re}\left[\delta_1+\delta_3\right]+\operatorname{Im}\left[\delta_6-\lambda\delta_1\right]\right\}+(e^{\pi\lambda}+1)\left\{\operatorname{Re}\left[\delta_2+\delta_4\right]+\operatorname{Im}\left[\delta_5-\lambda\delta_2\right]\right\},\vspace{0.3cm}\\
%&&+(e^{4\pi\lambda}+1)\operatorname{Im}\left[\delta_5\right]+(e^{\pi\lambda}-1)\operatorname{Im}\left[\delta_6\right]\vspace{0.3cm}
\end{array}$$
where
$$\begin{array}{rcl}
\delta_1&=&\dfrac{1}{8}A\left(\dfrac{\overline{\gamma_2}}{\lambda}+\dfrac{\gamma_1}{i+\lambda}+\dfrac{2\overline{\gamma_1}}{-i+2\lambda}+\dfrac{2\gamma_2}{i+2\lambda}\right);\vspace{0.3cm}\\
\delta_2&=&\dfrac{1}{8}A\left(\dfrac{4\gamma_1}{i+\lambda}+\dfrac{4\gamma_2}{i+\lambda}+\dfrac{4\overline{\gamma_1}}{i+\lambda}+\dfrac{4\overline{\gamma_2}}{i+\lambda}\right);\vspace{0.3cm}\\
\delta_3&=&\dfrac{A}{8(i+\lambda)(i+2\lambda)}(-\gamma_3+\overline{\gamma_4}+3i\gamma_3\lambda+2i\gamma_4\lambda-2i\overline{\gamma_3}\lambda-3i\overline{\gamma_4}\lambda\vspace{0.3cm}\\
&&+2\gamma_3\lambda^2+2\gamma_4\lambda^2-2\overline{\gamma_3}\lambda^2-2\overline{\gamma_4}\lambda^2);\vspace{0.3cm}\\
\delta_4&=&\dfrac{A}{8(i+\lambda)(i+2\lambda)}(4i\gamma_3\lambda+4i\gamma_4\lambda-4i\overline{\gamma_3}\lambda-4i\overline{\gamma_4}\lambda\vspace{0.3cm}\\
&&+8\gamma_3\lambda^2+8\gamma_4\lambda^2-8\overline{\gamma_3}\lambda^2-8\overline{\gamma_4}\lambda^2);\vspace{0.3cm}\\
\delta_5&=&\dfrac{1}{8}\left(\dfrac{4\overline{A}\gamma_3}{-i+\lambda}+\dfrac{4\overline{A}\gamma_4}{-i+\lambda}+\dfrac{4A\gamma_3}{i+\lambda}+\dfrac{4A\gamma_4}{i+\lambda}\right);\vspace{0.3cm}\\
\delta_6&=&\dfrac{1}{8}\left(\dfrac{A\gamma_3}{\lambda}+\dfrac{\overline{A}\gamma_4}{\lambda}+\dfrac{2\overline{A}\gamma_3}{-i+2\lambda}+\dfrac{2A\gamma_4}{i+2\lambda}\right).\vspace{0.3cm}\\
\end{array}$$
$$\begin{array}{rcl}
\widetilde{(\widetilde{R_2})^2R_3}(\pi)&=& \dfrac{1}{4}\displaystyle\operatorname{Re}\int_{0}^{\pi}e^{\lambda\theta}(\kappa+\overline{\kappa})^2\left[S_3(\varphi)+\frac{\lambda}{2}(-S_2^2(\varphi)+S_2(\varphi)\overline{S_2}(\varphi))\right]d\varphi\vspace{0.3cm}\\
&&-\dfrac{1}{4}\displaystyle\operatorname{Im}\int_{0}^{\pi}e^{2\lambda\theta}(\kappa+\overline{\kappa})^2\left[\lambda S_3(\varphi)+\frac{1}{2}S_2^2(\varphi)\right]d\varphi\vspace{0.3cm}\\
&=&\dfrac{1}{4}\displaystyle\operatorname{Re}\int_{0}^{\pi}e^{2\lambda\theta}(\kappa+\overline{\kappa})^2(1+\lambda i)S_3(\varphi)d\varphi\vspace{0.3cm}\\
&&+\dfrac{\lambda}{8}\displaystyle\operatorname{Re}\int_{0}^{\pi}e^{2\lambda\theta}(\kappa+\overline{\kappa})^2S_2(\varphi)\overline{S_2}(\varphi)d\varphi\vspace{0.3cm}\\
&&-\dfrac{1}{8}\displaystyle\operatorname{Im}\int_{0}^{\pi}e^{2\lambda\theta}(\kappa+\overline{\kappa})^2(1+\lambda i)(S_2(\varphi))^2d\varphi\vspace{0.3cm}\\
&=& \dfrac{1}{4}(e^{2\pi\lambda}-1)\left\{\operatorname{Re}\left[B\delta_7+\dfrac{\lambda}{2}\delta_{10}\right]-\dfrac{1}{2}\operatorname{Im}\left[A^2\delta_7\right]\right\}\vspace{0.3cm}\\
&&-\dfrac{1}{4}(-e^{\pi\lambda}-1)\left\{\operatorname{Re}\left[B\delta_8+\dfrac{\lambda}{2}\delta_9\right]-\dfrac{1}{2}\operatorname{Im}\left[A^2\delta_8\right]\right\},\vspace{0.3cm}\\
\end{array}$$
where
$$\begin{array}{rcl}
\delta_7&=&-\dfrac{i}{4(-i+\lambda)(i+\lambda)^2}\left(-2i\overline{A}^2-6\overline{A}^2\lambda+6i\overline{A}^2\lambda^2+2\overline{A}^2\lambda^3-\dfrac{4A\overline{A}}{2i+\lambda}\right.\vspace{0.2cm}\\
&&-\left.\dfrac{8A\overline{A}\lambda^2}{2i+\lambda}-\dfrac{4A\overline{A}\lambda^4}{2i+\lambda}+\dfrac{2A^2}{3i+\lambda}+\dfrac{8iA^2\lambda}{3i+\lambda}-\dfrac{12A^2\lambda^2}{3i+\lambda}-\dfrac{8iA^2\lambda^3}{3i+\lambda}+\dfrac{2A^2\lambda^4}{3i+\lambda}\right);\vspace{0.2cm}\\
\delta_8&=&-\dfrac{i}{4(-i+\lambda)(i+\lambda)^2}\left(-\dfrac{8A\overline{A}}{3i+\lambda}+\dfrac{8\overline{A}^2}{3i+\lambda}-\dfrac{32i\overline{A}^2\lambda}{3i+\lambda}-\dfrac{16A\overline{A}\lambda^2}{3i+\lambda}\right.\vspace{0.2cm}\\
&&-\dfrac{48\overline{A}^2\lambda^2}{3i+\lambda}+\dfrac{32i\overline{A}^2\lambda^3}{3i+\lambda}-\dfrac{8A\overline{A}\lambda^4}{3i+\lambda}+\dfrac{8\overline{A}^2\lambda^4}{3i+\lambda}+\dfrac{8A^2}{5i+\lambda}-\dfrac{8A\overline{A}}{5i+\lambda}+\dfrac{32iA^2\lambda}{5i+\lambda}\vspace{0.2cm}\\
&&-\left.\dfrac{48A^2\lambda^2}{5i+\lambda}-\dfrac{16A\overline{A}\lambda^2}{5i+\lambda}-\dfrac{32iA^2\lambda^3}{5i+\lambda}+\dfrac{8A^2\lambda^4}{5i+\lambda}-\dfrac{8A\overline{A}\lambda^4}{5i+\lambda}\right);\vspace{0.2cm}\\
\delta_9&=&-\dfrac{1}{2(i+\lambda)^2}A\overline{A}\left(\dfrac{4iA\overline{A}}{(-i+\lambda)^2}-\dfrac{4i\overline{A}^2}{(-i+\lambda)^2}-\dfrac{4A\overline{A}\lambda}{(-i+\lambda)^2}-\dfrac{12\overline{A}^2\lambda}{(-i+\lambda)^2}\right.\vspace{0.2cm}\\
&&+\dfrac{4iA\overline{A}\lambda^2}{(-i+\lambda)^2}+\dfrac{12i\overline{A}^2\lambda^2}{(-i+\lambda)^2}-\dfrac{4A\overline{A}\lambda^3}{(-i+\lambda)^2}+\dfrac{4\overline{A}^2\lambda^3}{(-i+\lambda)^2}-\dfrac{4A^2}{3i+\lambda}+\dfrac{4A\overline{A}}{3i+\lambda}\vspace{0.2cm}\\
&&-\left.\dfrac{8iA^2\lambda}{3i+\lambda}-\dfrac{8iA\overline{A}\lambda}{3i+\lambda}+\dfrac{4A^2\lambda^2}{3i+\lambda}-\dfrac{4A\overline{A}\lambda^2}{3i+\lambda}\right)\vspace{0.2cm};\\
\delta_{10}&=&-\dfrac{1}{2(i+\lambda)^2}A\overline{A}\left(-2iA\overline{A}-2A\overline{A}\lambda-\dfrac{A^2}{2i+\lambda}-\dfrac{2iA^2\lambda}{2i+\lambda}+\dfrac{A^2\lambda^2}{2i+\lambda}+\dfrac{6i\overline{A}^2}{(1+\lambda^2)^2}\right.\vspace{0.2cm}\\
&&-\left.\dfrac{\overline{A}^2}{\lambda(1+\lambda^2)^2}+\dfrac{15\overline{A}^2\lambda}{(1+\lambda^2)^2}-\dfrac{20i\overline{A}^2\lambda^2}{(1+\lambda^2)^2}-\dfrac{15\overline{A}^2\lambda^3}{(1+\lambda^2)^2}+\dfrac{6i\overline{A}^2\lambda^4}{(1+\lambda^2)^2}+\dfrac{\overline{A}^2\lambda^5}{(1+\lambda^2)^2}\right).
\end{array}$$
Thus,
$$\begin{array}{rcl}
\omega_5(\pi)&=&e^{\lambda\pi}u_5(\pi)\vspace{0.1cm}\\
&=&e^{\lambda\pi}\left(-\dfrac{5}{2}(u_2(\pi))^4+2(u_2(\pi))^2u_3(\pi)+\dfrac{3}{2}(u_3(\pi))^2+\widetilde{(\widetilde{R}_2)^2R_3}\right.\vspace{0.1cm}\\
&&-\left.2\widetilde{R}_2\widetilde{\widetilde{R}_3R_2}+4\widetilde{R}_2\widetilde{R}_4-2\widetilde{\widetilde{R}_4R_2}+\widetilde{R}_5\right)\\
&=&-\dfrac{5}{2}e^{-3\lambda\pi}(\omega_2(\pi))^4+2e^{-2\lambda\pi}(\omega_2(\pi))^2\omega_3(\pi)+\dfrac{3}{2}e^{-\lambda\pi}(\omega_3(\pi))^2\vspace{0.3cm}\\
 &&+e^{\lambda\pi}(e^{4\pi\lambda}-1)\left\{\operatorname{Re}\left[\xi_1+\xi_3+\xi_5+\xi_8+\xi_9+\xi_{10}+\xi_{11}+\xi_{12}-2(\delta_1+\delta_3)\right]\right.\vspace{0.3cm}\\
&&+\left.\operatorname{Im}\left[\xi_6+\xi_7-2(\delta_6-\lambda\delta_1
)\right]\right\}+e^{\lambda\pi}(-e^{4\pi\lambda}-1)\left\{\operatorname{Re}\left[\xi_2\right]-\operatorname{Im}\left[\xi_4\right]\right\}\vspace{0.3cm}\\   
 &&-\omega_2(\pi)(-e^{3\pi\lambda}-1)\left\{\operatorname{Re}\left[\eta_1-\lambda\eta_5-4(\gamma_1+\gamma_2)\right]-\operatorname{Im}\left[\lambda\eta_1+\eta_3+4(\gamma_3+\gamma_4)\right]\right\}\vspace{0.3cm}\\
&& +\dfrac{1}{4}e^{\lambda\pi}(e^{2\pi\lambda}-1)\left\{\operatorname{Re}\left[B\delta_7+\dfrac{\lambda}{2}\delta_{10}\right]-\dfrac{1}{2}\operatorname{Im}\left[A^2\delta_7\right]\right\}\vspace{0.3cm}\\
&&-\dfrac{1}{4}e^{\lambda\pi}(-e^{\pi\lambda}-1)\left\{\operatorname{Re}\left[B\delta_8+\dfrac{\lambda}{2}\delta_9-8(\delta_2+\delta_4)\right]-\dfrac{1}{2}\operatorname{Im}\left[A^2\delta_8+16(\delta_5-\lambda\delta_2)\right]\right\}\vspace{0.3cm}\\
&&+\omega_2(\pi)(-e^{\pi\lambda}-1)\left\{\operatorname{Re}\left[\eta_2-\lambda\eta_6\right]-\operatorname{Im}\left[\lambda\eta_2+\eta_4\right]\right\}.\\
\end{array}$$

\section{Acknowledgements}
This article was possible thanks to the scholarship granted from the Brazilian Federal Agency for Support and Evaluation of Graduate Education (CAPES), in the scope of the Program CAPES-Print, process number 88887.310463/2018-00, International Cooperation Project number 88881.310741/2018-01. 

Armengol Gasull is partially supported by Spanish State Research Agency, through the projects PID2019-104658GB-I00 and PID2022-136613NB-I00 grants and the Severo Ochoa and María de Maeztu Program for Centers and Units of Excellence in R$\&$D (CEX2020-001084-M), and grant 2021-SGR-00113 from AGAUR. Generalitat de Catalunya.

Gabriel Rondón is supported by São Paulo Research Foundation (FAPESP) grants 2020/06708-9 and 2022/12123-9. 

 Paulo Ricardo da Silva is also partially supported by São Paulo Research Foundation (FAPESP) grant 2019/10269-3.
%====================================================
%====================================================
%====================================================
%\newpage
%\appendix

%====================================================
%====================================================
%====================================================

%\begin{thebibliography}{99}

%\bibitem{SilvaSarmientoNovaes} P.R. da Silva, I.S. Meza-Sarmiento, D.D. Novaes. \emph{Nonlinear Sliding of Discontinuous Vector Fields and Singular Perturbation}. \textbf{Differ Equ Dyn Syst} (2018), DOI https://doi.org/10.1007/s12591-018-0439-1.

%\bibitem{SilvaMoraes} P.R. da Silva, J.R. de Moraes. \emph{Piecewise-Smooth Slow–Fast Systems}. \textbf{J. Dyn. Control Syst}. 27 (2021), pp 67--85.

%\bibitem{SotoTeixeira} J. Sotomayor, M.A. Teixeira. \emph{Regularization of discontinuous vector fields}. International Conference on Differential Equations, Lisboa, 1996, pp 207--223.

%\bibitem{Szmolyan} P. Szmolyan. \emph{Transversal heteroclinic and homoclinic orbits in singular perturbation problems}. \textbf{J. Diff. Eq.} 92 (1991), pp 252--281.

%\bibitem{Wiggins} S. Wiggins. \emph{Normally Hyperbolic Invariant Manifolds in Dynamical Systems}. Springer, 1994.	

%\end{thebibliography}

\bibliographystyle{abbrv}
\bibliography{references1}

\end{document}